\documentclass[reqno]{amsart}
\usepackage{amssymb,latexsym}
\usepackage{graphicx}
\usepackage{fancyhdr}
\numberwithin{equation}{section}
\newtheorem{thm}{Theorem}[section]
\newtheorem{theorem}[thm]{Theorem}
\newtheorem{remark}[thm]{Remark}
\newtheorem{lemma}[thm]{Lemma}

\newtheorem{corollary}[thm]{Corollary}

\begin{document}

\setcounter{page}{1}

\title[Reciprocity formulas]{Reciprocity formulas for certain generalized Hardy-Berndt sums}
\author{Yuan He}
\address{School of Mathematics and Big Data, Neijiang Normal University, Neijiang 641100, Sichuan, People's Republic of China}
\email{hyyhe@aliyun.com}

\begin{abstract}
In this paper, we introduce the generalized Hardy-Berndt sums. We establish some formulas of products of the Bernoulli and Euler functions by using the Fourier series technique and some properties of the periodic-zeta and Lerch-zeta functions. As applications of the results, we give some reciprocity formulas for the generalized Hardy-Berndt sums. It turns out that Hardy's reciprocity formulas follow as special cases and that Goldberg's three-term reciprocity formulas are recovered.
\end{abstract}

\maketitle

\section{Introduction}

Let $\mathbb{N}$ be the set of positive integers, $\mathbb{N}_{0}$ the set of non-negative integers, $\mathbb{Z}$ the set of integers, $\mathbb{R}$ the set of real numbers, and $\mathbb{C}$ the set of complex numbers. Denote by $(a,b)$ the greatest common divisor of $a,b\in\mathbb{Z}$, by $\{x\}$ the fractional part of $x\in\mathbb{R}$, by $[x]$ the greatest integer function (also called the floor function), i.e., $[x]=x-\{x\}$ for $x\in\mathbb{R}$, and by $((x))$ the sawtooth function (also called the first Bernoulli function) given by
\begin{equation*}
((x))=\begin{cases}
\{x\}-\frac{1}{2},  &\text{if $x\in\mathbb{R}\setminus\mathbb{Z}$},\\
0,  &\text{if $x\in\mathbb{Z}$}.
\end{cases}
\end{equation*}
The classical Dedekind sums $s(a,b)$, arising in the transformation formulas of the logarithm of the Dedekind eta function, is defined for $a\in\mathbb{Z}$, $b\in\mathbb{N}$ by
\begin{equation}\label{eq1.1}
s(a,b)=\sum_{r=1}^{b}\biggl(\biggl(\frac{r}{b}\biggl)\biggl)\biggl(\biggl(\frac{ar}{b}\biggl)\biggl).
\end{equation}
It is well known that one of the most intriguing and important features of the sums \eqref{eq1.1} is the following Dedekind's \cite{dedekind} reciprocity theorem. If $a,b\in\mathbb{N}$ are such that $(a,b)=1$, then
\begin{equation}\label{eq1.2}
s(a,b)+s(b,a)=-\frac{1}{4}+\frac{1}{12}\biggl(\frac{a}{b}+\frac{b}{a}+\frac{1}{ab}\biggl).
\end{equation}

In 1905, Hardy \cite{hardy} first provided a different proof of the reciprocity formula \eqref{eq1.2} which does not depend on the theory of the Dedekind eta function, and stated without proofs some reciprocity formulas for arithmetic sums similar to \eqref{eq1.1}. For example, Hardy \cite[pp. 122--123]{hardy} showed that for $a,b\in\mathbb{N}$ with $(a,b)=1$, if $2\nmid(a+b)$ then
\begin{equation}\label{eq1.3}
S(a,b)+S(b,a)=1;
\end{equation}
if $2\mid a$ then
\begin{equation}\label{eq1.4}
s_{1}(a,b)-2s_{2}(b,a)=\frac{1}{2}-\frac{1}{2}\biggl(\frac{1}{ab}+\frac{b}{a}\biggl);
\end{equation}
if $2\nmid b$ then
\begin{equation}\label{eq1.5}
2s_{3}(a,b)-s_{4}(b,a)=1-\frac{a}{b};
\end{equation}
if $2\nmid a$ and $2\nmid b$ then
\begin{equation}\label{eq1.6}
s_{5}(a,b)+s_{5}(b,a)=\frac{1}{2}-\frac{1}{2ab},
\end{equation}
where $S(a,b)$, $s_{1}(a,b)$, $s_{2}(a,b)$, $s_{3}(a,b)$, $s_{4}(a,b)$, $s_{5}(a,b)$ are the Hardy sums (also known as the Hardy-Berndt sums) given for $a\in\mathbb{Z}$, $b\in\mathbb{N}$ by
\begin{equation*}
S(a,b)=\sum_{r=1}^{b-1}(-1)^{r+1+[\frac{ar}{b}]},\quad
s_{1}(a,b)=\sum_{r=1}^{b}(-1)^{[\frac{ar}{b}]}\biggl(\biggl(\frac{r}{b}\biggl)\biggl),
\end{equation*}
\begin{equation*}
s_{2}(a,b)=\sum_{r=1}^{b}(-1)^{r}\biggl(\biggl(\frac{ar}{b}\biggl)\biggl)\biggl(\biggl(\frac{r}{b}\biggl)\biggl),\quad
s_{3}(a,b)=\sum_{r=1}^{b}(-1)^{r}\biggl(\biggl(\frac{ar}{b}\biggl)\biggl),
\end{equation*}
\begin{equation*}
s_{4}(a,b)=\sum_{r=1}^{b-1}(-1)^{[\frac{ar}{b}]},\quad
s_{5}(a,b)=\sum_{r=1}^{b}(-1)^{r+[\frac{ar}{b}]}\biggl(\biggl(\frac{r}{b}\biggl)\biggl).
\end{equation*}
Just as Dedekind \cite{dedekind} deduced the reciprocity formula \eqref{eq1.2} from his transformation formula of the logarithm of the Dedekind eta function, Berndt \cite{berndt2} in 1978 used his transformation formulas of the logarithms of the classical theta-functions to demonstrate Hardy's reciprocity formulas \eqref{eq1.3}, \eqref{eq1.4} and \eqref{eq1.5}. In 1981, Goldberg \cite{goldberg} used Berndt's \cite{berndt2} transformation formulas to demonstrate Hardy's reciprocity formula \eqref{eq1.6}. In fact, Goldberg \cite{goldberg} also presented some three-term reciprocity formulas for the above Hardy-Berndt sums. After that, Apostol and Vu \cite{apostol2}, Berndt and Dieter \cite{berndt3}, Berndt and Goldberg \cite{berndt4}, Dieter \cite{dieter}, Sitaramachandrarao \cite{sitaramachandrarao}, et al. provided some different proofs of Hardy's reciprocity formulas \eqref{eq1.3}--\eqref{eq1.6}. In particular, Berndt and Goldberg \cite{berndt4} gave the infinite and finite series representations of the Hardy-Berndt sums, and Sitaramachandrarao \cite{sitaramachandrarao} expressed each one of the Hardy-Berndt sums explicitly in terms of the classical Dedekind sums \eqref{eq1.1}. Pettet and Sitaramachandrarao \cite{pettet} and Simsek \cite{simsek} also adopted different strategies to reobtain Goldberg's \cite{goldberg} three-term reciprocity formulas.

The object of the present paper is to give some new extensions of Hardy's reciprocity formulas \eqref{eq1.3}--\eqref{eq1.6}. To do so, following the work of Hall, Wilson and Zagier \cite{hall}, we study the following generalized Hardy-Berndt sums:
\begin{equation}\label{eq1.7}
S_{m,n}\left(
\begin{matrix}
a & b & c \\
x & y & z
\end{matrix}
\right)=\sum_{r=1}^{c}(-1)^{r}\overline{E}_{m}\biggl(\frac{a(r+z)}{c}-x\biggl)\overline{E}_{n}\biggl(\frac{b(r+z)}{c}-y\biggl),
\end{equation}
\begin{equation}\label{eq1.8}
s_{m,n}^{(1)}\left(
\begin{matrix}
a & b & c \\
x & y & z
\end{matrix}
\right)=\sum_{r=1}^{c}\overline{E}_{m}\biggl(\frac{a(r+z)}{c}-x\biggl)\overline{B}_{n}\biggl(\frac{b(r+z)}{c}-y\biggl),
\end{equation}
\begin{equation}\label{eq1.9}
s_{m,n}^{(2)}\left(
\begin{matrix}
a & b & c \\
x & y & z
\end{matrix}
\right)=\sum_{r=1}^{c}(-1)^{r}\overline{B}_{m}\biggl(\frac{a(r+z)}{c}-x\biggl)\overline{B}_{n}\biggl(\frac{b(r+z)}{c}-y\biggl),
\end{equation}
\begin{equation}\label{eq1.10}
s_{m,n}^{(4)}\left(
\begin{matrix}
a & b & c \\
x & y & z
\end{matrix}
\right)=\sum_{r=1}^{c}\overline{E}_{m}\biggl(\frac{a(r+z)}{c}-x\biggl)\overline{E}_{n}\biggl(\frac{b(r+z)}{c}-y\biggl),
\end{equation}
and
\begin{equation}\label{eq1.11}
s_{m,n}^{(5)}\left(
\begin{matrix}
a & b & c \\
x & y & z
\end{matrix}
\right)=\sum_{r=1}^{c}(-1)^{r}\overline{E}_{m}\biggl(\frac{a(r+z)}{c}-x\biggl)\overline{B}_{n}\biggl(\frac{b(r+z)}{c}-y\biggl),
\end{equation}
where $m,n\in\mathbb{N}_{0}$, $a,b\in\mathbb{Z}$, $c\in\mathbb{N}$, $x,y,z\in\mathbb{R}$, $\overline{B}_{n}(x)$ is the $n$-th Bernoulli function given for $n\in\mathbb{N}_{0}$, $x\in\mathbb{R}$ by
\begin{equation}\label{eq1.12}
\overline{B}_{0}(x)=1, \quad\overline{B}_{1}(x)=((x)),\quad\overline{B}_{n}(x)=B_{n}(\{x\})\quad(n\geq2),
\end{equation}
$\overline{E}_{n}(x)$ is the $n$-th Euler function given for $n\in\mathbb{N}_{0}$, $x\in\mathbb{R}$ by
\begin{equation}\label{eq1.13}
\overline{E}_{0}(x)=\begin{cases}
(-1)^{[x]},  &\text{if $x\in\mathbb{R}\setminus\mathbb{Z}$},\\
0,  &\text{if $x\in\mathbb{Z}$},
\end{cases} \quad\overline{E}_{n}(x)=(-1)^{[x]}E_{n}(\{x\})\quad(n\geq1),
\end{equation}
in which $B_{n}(x)$ is the Bernoulli polynomial of degree $n$, and $E_{n}(x)$ is the Euler polynomial of degree $n$. By making use of the Fourier series technique and some properties of the periodic-zeta function and the Lerch-zeta function, we establish some formulas of products of the Bernoulli function and the Euler function (see Theorems \ref{thm2.1}, \ref{thm2.5}, \ref{thm2.8} and \ref{thm2.10} below). As applications of the results, we obtain some reciprocity formulas for the generalized Hardy-Berndt sums. We also show that Hardy's reciprocity formulas \eqref{eq1.3}--\eqref{eq1.6} follow as special cases, and Goldberg's \cite{goldberg} three-term reciprocity formulas are recovered.

This paper is organized as follows. In the second section, we state our main results and their applications. In the third section, we give some auxiliary lemmas. The fourth section is devoted to the detailed proofs of Theorems \ref{thm2.1}, \ref{thm2.5}, \ref{thm2.8} and \ref{thm2.10}.

\section{Statement of main results}

For convenience, in the following we always denote by $\mathrm{i}$ the square root of $-1$ such that $\mathrm{i}^{2}=-1$, and $\Gamma(s)$ the gamma function defined on $s\in\mathbb{C}$. We also write, for $l,k\in\mathbb{N}_{0}$, $\delta_{l,k}$ as the Kronecker delta function given by $\delta_{l,k}=1$ or $0$ according to $l=k$ or $l\not=k$; for $x\in\mathbb{R}\setminus\{0\}$, $\mathrm{sgn}(x)$ as the sign of $x$ given by $\mathrm{sgn}(x)=x/|x|$. For the sake of convergence, we denote that
\begin{equation*}
\sum_{k=-\infty}^{+\infty}\frac{1}{k+a}=\lim_{N\rightarrow\infty}\sum_{k=-N}^{N}\frac{1}{k+a}\quad(a\not\in\mathbb{Z}).
\end{equation*}
This ensures that the $n$-th Bernoulli function $\overline{B}_{n}(x)$ defined in \eqref{eq1.12} can be given for $n\in\mathbb{N}$, $x\in\mathbb{R}$ by the Fourier series (see, e.g., \cite[Corollary 3.2]{he3} or \cite[Corollary 3.3]{he3})
\begin{equation}\label{eq2.1}
\overline{B}_{n}(x)=-\frac{n!}{(2\pi\mathrm{i})^{n}}\sideset{}{'}\sum_{k=-\infty}^{+\infty}\frac{e^{2\pi\mathrm{i}kx}}{k^{n}},
\end{equation}
the $n$-th Euler function $\overline{E}_{n}(x)$ defined in \eqref{eq1.13} can be given for $n\in\mathbb{N}_{0}$, $x\in\mathbb{R}$ by the Fourier series (see, e.g., \cite[Corollary 3.2]{he3} or \cite[Corollary 3.3]{he3})
\begin{equation}\label{eq2.2}
\overline{E}_{n}(x)=\frac{2n!}{(2\pi\mathrm{i})^{n+1}}\sum_{k=-\infty}^{+\infty}\frac{e^{2\pi\mathrm{i}(k-\frac{1}{2})x}}{\bigl(k-\frac{1}{2}\bigl)^{n+1}}.
\end{equation}
Here the dash appearing in \eqref{eq2.1} denotes throughout that undefined terms are excluded from the sum. Our main results are as follows.

\begin{theorem}\label{thm2.1} Let $m\in\mathbb{N}_{0}$, $n\in\mathbb{N}$, $a,b\in\mathbb{Z}\setminus\{0\}$ with $2\mid a$ and $2\nmid b$. Then, for $x,y,z\in\mathbb{R}$,
\begin{eqnarray}\label{eq2.3}
&&\overline{E}_{m}(ax+y)\overline{B}_{n}(bx+z)\nonumber\\
&&=-nb^{n-1}\mathrm{sgn}(b)
\sum_{j=0}^{m}\binom{m}{j}\frac{(-1)^{j}a^{m-j}}{m+n-j}\sum_{l=1}^{|b|}\overline{E}_{j}\biggl(\frac{a(l+z)}{b}-y\biggl)\nonumber\\
&&\qquad\times\overline{B}_{m+n-j}\biggl(x+\frac{l+z}{b}\biggl)\nonumber\\
&&\quad-2a^{m}\mathrm{sgn}(a)\sum_{j=1}^{n}\binom{n}{j}\frac{(-1)^{j}b^{n-j}}{m+n+1-j}\sum_{l=1}^{|a|}(-1)^{l}\overline{B}_{j}\biggl(\frac{b(l+y)}{a}-z\biggl)\nonumber\\
&&\qquad\times\overline{B}_{m+n+1-j}\biggl(x+\frac{l+y}{a}\biggl)\nonumber\\
&&\quad+\frac{(-1)^{n+1}m!n!(a,b)^{m+n+1}}{a^{n}b^{m+1}(m+n)!}\overline{E}_{m+n}\biggl(\frac{by-az}{(a,b)}\biggl)+\frac{b^{n}}{a^{n}}\overline{E}_{m+n}(ax+y)\nonumber\\
&&\quad+\frac{1}{2}\delta_{0,m}\delta_{1,n}\mathrm{sgn}(ab)(-1)^{ax+y}\delta_{\mathbb{Z}}(ax+y)\delta_{\mathbb{Z}}(bx+z),
\end{eqnarray}
where $\delta_{\mathbb{Z}}(x)=1$ or $0$ according to $x\in\mathbb{Z}$ or $x\not\in\mathbb{Z}$.
\end{theorem}

If we take $m=x=y=z=0$ and $n=1$ in Theorem \ref{thm2.1}, then by $\overline{E}_{1}(0)=-1/2$, we obtain the reciprocity formula \eqref{eq1.4}.

We next use Theorem \ref{thm2.1} to give a more general reciprocity formula among the generalized Hardy-Berndt sums \eqref{eq1.8} and \eqref{eq1.9}.

\begin{corollary}\label{cor2.2} Let $m\in\mathbb{N}_{0}$, $n,a,b,c\in\mathbb{N}$ with $2\mid a$ and $2\nmid b$. Then, for $x,y,z\in\mathbb{R}$,
\begin{eqnarray}\label{eq2.4}
&&nb^{n-1}\sum_{j=0}^{m}\binom{m}{j}\frac{(-1)^{j}a^{m-j}}{c^{m+n-1-j}(m+n-j)}s_{j,m+n-j}^{(1)}\left(
\begin{matrix}
a & c & b \\
y & x & z
\end{matrix}
\right)\nonumber\\
&&\quad+2a^{m}\sum_{j=1}^{n}\binom{n}{j}\frac{(-1)^{j}b^{n-j}}{c^{m+n-j}(m+n+1-j)}s_{j,m+n+1-j}^{(2)}\left(
\begin{matrix}
b & c & a \\
z & x & y
\end{matrix}
\right)\nonumber\\
&&=-s_{m,n}^{(1)}\left(
\begin{matrix}
-a & -b & c \\
-y & -z & x
\end{matrix}
\right)+\frac{b^{n}}{a^{n}}s_{m+n,0}^{(1)}\left(
\begin{matrix}
-a & -b & c \\
-y & -z & x
\end{matrix}
\right)\nonumber\\
&&\quad+\frac{(-1)^{n+1}m!n!(a,b)^{m+n+1}c}{a^{n}b^{m+1}(m+n)!}\overline{E}_{m+n}\biggl(\frac{by-az}{(a,b)}\biggl)\nonumber\\
&&\quad+\frac{1}{2}\delta_{0,m}\delta_{1,n}\sum_{r=1}^{c}(-1)^{\frac{a(r+x)}{c}-y}
\delta_{\mathbb{Z}}\biggl(\frac{a(r+x)}{c}-y\biggl)\delta_{\mathbb{Z}}\biggl(\frac{b(r+x)}{c}-z\biggl).
\end{eqnarray}
\end{corollary}

\begin{proof}
Since the Bernoulli function satisfies the so-called ``distribution property'' (Raabe's formula; see \cite{raabe} or \cite[Equation (3.43)]{he3})
\begin{equation}\label{eq2.5}
\overline{B}_{n}(qx)=q^{n-1}\sum_{r=0}^{q-1}\overline{B}_{n}\biggl(x+\frac{r}{q}\biggl),
\end{equation}
where $q,n\in\mathbb{N}$, $x\in\mathbb{R}$, we see from \eqref{eq1.12} and \eqref{eq2.5} that for $n,c\in\mathbb{N}$, $x\in\mathbb{R}$,
\begin{equation}\label{eq2.6}
\sum_{r=1}^{c}\overline{B}_{n}\biggl(\frac{-r+x}{c}\biggl)=\sum_{r=0}^{c-1}\overline{B}_{n}\biggl(\frac{-(c-r)+x}{c}\biggl)=c^{1-n}\overline{B}_{n}(x).
\end{equation}
Hence, replacing $x$ by $-(r+x)/c$ and then applying the operation $\sum_{r=1}^{c}$ to both sides of \eqref{eq2.3},  we get from \eqref{eq2.6} that
\begin{eqnarray*}
&&\sum_{r=1}^{c}\overline{E}_{m}\biggl(-\frac{a(r+x)}{c}+y\biggl)\overline{B}_{n}\biggl(-\frac{b(r+x)}{c}+z\biggl)\\
&&=-nb^{n-1}\sum_{j=0}^{m}\binom{m}{j}\frac{(-1)^{j}a^{m-j}}{c^{m+n-1-j}(m+n-j)}
\sum_{l=1}^{b}\overline{E}_{j}\biggl(\frac{a(l+z)}{b}-y\biggl)\\
&&\qquad\times\overline{B}_{m+n-j}\biggl(\frac{c(l+z)}{b}-x\biggl)\\
&&\quad-2a^{m}\sum_{j=1}^{n}\binom{n}{j}\frac{(-1)^{j}b^{n-j}}{c^{m+n-j}(m+n+1-j)}\sum_{l=1}^{a}(-1)^{l}\overline{B}_{j}\biggl(\frac{b(l+y)}{a}-z\biggl)\\
&&\qquad\times\overline{B}_{m+n+1-j}\biggl(\frac{c(l+y)}{a}-x\biggl)\\
&&\quad+\frac{(-1)^{n+1}m!n!(a,b)^{m+n+1}c}{a^{n}b^{m+1}(m+n)!}\overline{E}_{m+n}\biggl(\frac{by-az}{(a,b)}\biggl)\\
&&\quad+\frac{b^{n}}{a^{n}}\sum_{r=1}^{c}\overline{E}_{m+n}\biggl(-\frac{a(r+x)}{c}+y\biggl)\\
&&\quad+\frac{1}{2}\delta_{0,m}\delta_{1,n}\sum_{r=1}^{c}(-1)^{-\frac{a(r+x)}{c}+y}
\delta_{\mathbb{Z}}\biggl(-\frac{a(r+x)}{c}+y\biggl)\delta_{\mathbb{Z}}\biggl(-\frac{b(r+x)}{c}+z\biggl),
\end{eqnarray*}
as desired.
\end{proof}

As a special case of Corollary \ref{cor2.2}, we obtain the following extension of the reciprocity formula \eqref{eq1.4}.

\begin{corollary}\label{cor2.3} Let $a,b,c\in\mathbb{N}$ with $2\mid a$, $2\nmid b$ and $(a,b)=(a,c)=1$. Then
\begin{equation}\label{eq2.7}
s_{1}(a,b,c)-2s_{2}(b,c,a)+s_{1}(a,c,b)=\frac{1}{2}-\frac{1}{2}\biggl(\frac{c}{ab}+\frac{b}{ac}\biggl),
\end{equation}
where $s_{1}(a,b,c)$, $s_{2}(a,b,c)$ are given for $a,b\in\mathbb{Z}$, $c\in\mathbb{N}$ by
\begin{equation*}
s_{1}(a,b,c)=\sum_{r=1}^{c}(-1)^{[\frac{ar}{c}]}\biggl(\biggl(\frac{br}{c}\biggl)\biggl),
\end{equation*}
and
\begin{equation*}
s_{2}(a,b,c)=\sum_{r=1}^{c}(-1)^{r}\biggl(\biggl(\frac{ar}{c}\biggl)\biggl)\biggl(\biggl(\frac{br}{c}\biggl)\biggl).
\end{equation*}
\end{corollary}

\begin{proof}
It is easily seen from \eqref{eq2.1} that for $n\in\mathbb{N}$, $x\in\mathbb{R}$,
\begin{equation}\label{eq2.8}
\overline{B}_{n}(-x)=(-1)^{n}\overline{B}_{n}(x).
\end{equation}
Since $[-x]=-[x]-1$ in the case when $x\not\in\mathbb{Z}$, by taking $m=x=y=z=0$ and $n=1$ in Corollary \ref{cor2.2}, in light of \eqref{eq1.13} and \eqref{eq2.8}, we know that for $a,b,c\in\mathbb{N}$ with $2\mid a$, $2\nmid b$ and $(a,b)=(a,c)=1$,
\begin{eqnarray}\label{eq2.9}
&&s_{1}(a,b,c)-2s_{2}(b,c,a)+s_{1}(a,c,b)\nonumber\\
&&=\frac{1}{2}-\frac{c}{2ab}+\frac{b}{a}\sum_{r=1}^{c}\overline{E}_{1}\biggl(-\frac{ar}{c}\biggl)\nonumber\\
&&=\frac{1}{2}-\frac{c}{2ab}+\frac{b}{a}\sum_{r=0}^{c-1}\overline{E}_{1}\biggl(-\frac{a(c-r)}{c}\biggl)\nonumber\\
&&=\frac{1}{2}-\frac{c}{2ab}+\frac{b}{a}\sum_{r=0}^{c-1}\overline{E}_{1}\biggl(\frac{ar}{c}\biggl).
\end{eqnarray}
Note that from \eqref{eq2.2} and the familiar geometric sum stated in \cite[Theorem 8.1]{apostol1}, we discover that for $n\in\mathbb{N}_{0}$, $a\in\mathbb{Z}\setminus\{0\}$, $c\in\mathbb{N}$, $x\in\mathbb{R}$ with $2\mid a$ and $(a,c)=1$,
\begin{eqnarray}\label{eq2.10}
\sum_{r=0}^{c-1}\overline{E}_{n}\biggl(x+\frac{ar}{c}\biggl)
&=&\frac{2n!}{(2\pi\mathrm{i})^{n+1}}\sum_{k=-\infty}^{+\infty}\frac{e^{2\pi\mathrm{i}(k-\frac{1}{2})x}}{\bigl(k-\frac{1}{2}\bigl)^{n+1}}
\sum_{r=0}^{c-1}e^{\frac{2\pi\mathrm{i}(2k-1)\frac{a}{2}r}{c}}\nonumber\\
&=&\frac{2cn!}{(2\pi\mathrm{i})^{n+1}}\sum_{\substack{k=-\infty\\c\mid(2k-1)\frac{a}{2}}}^{+\infty}\frac{e^{2\pi\mathrm{i}(k-\frac{1}{2})x}}{\bigl(k-\frac{1}{2}\bigl)^{n+1}}\nonumber\\
&=&\frac{2cn!}{(2\pi\mathrm{i})^{n+1}}\sum_{\substack{k=-\infty\\c\mid(2k-1)}}^{+\infty}\frac{e^{2\pi\mathrm{i}(k-\frac{1}{2})x}}{\bigl(k-\frac{1}{2}\bigl)^{n+1}}\nonumber\\
&=&\frac{2n!}{(2\pi\mathrm{i})^{n+1}c^{n}}\sum_{k=-\infty}^{+\infty}\frac{e^{2\pi\mathrm{i}(k-\frac{1}{2})cx}}{\bigl(k-\frac{1}{2}\bigl)^{n+1}}\nonumber\\
&=&c^{-n}\overline{E}_{n}(cx).
\end{eqnarray}
Therefore, by applying \eqref{eq2.10} to the right-hand side of \eqref{eq2.9}, the desired result follows immediately.
\end{proof}

\begin{remark}\label{rem2.4}
By a similar consideration as in \eqref{eq2.10}, one can easily see that the Euler function satisfies the so-called ``distribution property'' (see, e.g., \cite[Equation (3.44)]{he3})
\begin{equation}\label{eq2.11}
\overline{E}_{n}(qx)=q^{n}\sum_{r=0}^{q-1}(-1)^{r}\overline{E}_{n}\biggl(x+\frac{r}{q}\biggl),
\end{equation}
where $n\in\mathbb{N}_{0}$, $q\in\mathbb{N}$, $x\in\mathbb{R}$ with $2\nmid q$. A generalization of \eqref{eq2.11} can be also found in \eqref{eq2.18} or \eqref{eq4.16} below.
\end{remark}

It should be noted that replacing $a$ by $b$, $b$ by $c$ and $c$ by $a$ in Corollary \ref{cor2.3}, and then using Pettet and Sitaramachandrarao's results stated in \cite[Equations (3.10), (3.11) and (3.12]{pettet}, we rediscover Goldberg's \cite[Theorem 5.12]{goldberg} three-term reciprocity formula
\begin{equation}\label{eq2.12}
s_{1}(bc',a)-2s_{2}(ca',b)+s_{3}(ab',c)=\frac{1}{2}-\frac{1}{2}\biggl(\frac{a}{bc}+\frac{c}{ab}\biggl),
\end{equation}
where $a,b,c\in\mathbb{N}$ with $2\mid b$ and $(a,b)=(b,c)=(a,c)=1$, $a',b',c'\in\mathbb{Z}$ satisfy $aa'\equiv 1$ (mod $b$), $bb'\equiv 1$ (mod $c$), $cc'\equiv 1$ (mod $a$).

\begin{theorem}\label{thm2.5} Let $m\in\mathbb{N}_{0}$, $n\in\mathbb{N}$, $a,b\in\mathbb{Z}\setminus\{0\}$ with $2\nmid a$. Then, for $x,y,z\in\mathbb{R}$,
\begin{eqnarray}\label{eq2.13}
&&\overline{E}_{m}(ax+y)\overline{B}_{n}(bx+z)\nonumber\\
&&=\frac{nb^{n-1}\mathrm{sgn}(b)}{2}\sum_{j=0}^{m}
\binom{m}{j}(-1)^{j}a^{m-j}\sum_{l=1}^{|b|}\overline{E}_{j}\biggl(\frac{a(l+z)}{b}-y\biggl)\nonumber\\
&&\qquad\times\overline{E}_{m+n-1-j}\biggl(x+\frac{l+z}{b}\biggl)\nonumber\\
&&\quad+a^{m}\mathrm{sgn}(a)\sum_{j=1}^{n}\binom{n}{j}(-1)^{j}b^{n-j}\sum_{l=1}^{|a|}(-1)^{l}
\overline{B}_{j}\biggl(\frac{b(l+y)}{a}-z\biggl)\nonumber\\
&&\qquad\times\overline{E}_{m+n-j}\biggl(x+\frac{l+y}{a}\biggl)\nonumber\\
&&\quad+\frac{b^{n}}{a^{n}}\overline{E}_{m+n}(ax+y)\nonumber\\
&&\quad+\frac{1}{2}\delta_{0,m}\delta_{1,n}\mathrm{sgn}(ab)(-1)^{ax+y}\delta_{\mathbb{Z}}(ax+y)\delta_{\mathbb{Z}}(bx+z),
\end{eqnarray}
where $\delta_{\mathbb{Z}}(x)$ is as in \eqref{eq2.3}.
\end{theorem}

By taking $m=x=y=z=0$ and $n=1$ in Theorem \ref{thm2.5}, and then swapping $a$ and $b$, we get the reciprocity formula \eqref{eq1.5}.
More generally, we have the following reciprocity formula among the generalized Hardy-Berndt sums \eqref{eq1.10} and \eqref{eq1.11}.

\begin{corollary}\label{cor2.6} Let $m\in\mathbb{N}_{0}$, $n,a,b,c\in\mathbb{N}$ with $2\nmid a$ and $2\nmid c$. Then, for $x,y,z\in\mathbb{R}$,
\begin{eqnarray}\label{eq2.14}
&&\frac{nb^{n-1}}{2}\sum_{j=0}^{m}
\binom{m}{j}\frac{(-1)^{j}a^{m-j}}{c^{m+n-1-j}}s_{j,m+n-1-j}^{(4)}\left(
\begin{matrix}
a & c & b \\
y & x & z
\end{matrix}
\right)\nonumber\\
&&\quad+a^{m}\sum_{j=1}^{n}\binom{n}{j}\frac{(-1)^{j}b^{n-j}}{c^{m+n-j}}s_{m+n-j,j}^{(5)}\left(
\begin{matrix}
c & b & a \\
x & z & y
\end{matrix}
\right)\nonumber\\
&&=s_{m,n}^{(5)}\left(
\begin{matrix}
-a & -b & c \\
-y & -z & x
\end{matrix}
\right)-\frac{b^{n}}{a^{n}}s_{m+n,0}^{(5)}\left(
\begin{matrix}
-a & -b & c \\
-y & -z & x
\end{matrix}
\right)\nonumber\\
&&\quad-\frac{1}{2}\delta_{0,m}\delta_{1,n}\sum_{r=1}^{c}(-1)^{\frac{a(r+x)}{c}-y-r}\delta_{\mathbb{Z}}\biggl(\frac{a(r+x)}{c}-y\biggl)\nonumber\\
&&\qquad\times\delta_{\mathbb{Z}}\biggl(\frac{b(r+x)}{c}-z\biggl).
\end{eqnarray}
\end{corollary}

\begin{proof}
By \eqref{eq1.13} and \eqref{eq2.11}, we see that for $n\in\mathbb{N}_{0}$, $c\in\mathbb{N}$, $x\in\mathbb{R}$ with $2\nmid c$,
\begin{equation}\label{eq2.15}
\sum_{r=1}^{c}(-1)^{r}\overline{E}_{n}\biggl(\frac{-r+x}{c}\biggl)=\sum_{r=0}^{c-1}(-1)^{c-r}\overline{E}_{n}\biggl(\frac{-(c-r)+x}{c}\biggl)=c^{-n}\overline{E}_{n}(x).
\end{equation}
Thus, replacing $x$ by $-(r+x)/c$ and then applying the operation $\sum_{r=1}^{c}(-1)^{r}$ to both sides of \eqref{eq2.13}, we know from \eqref{eq2.15} that
\begin{eqnarray*}
&&\sum_{r=1}^{c}(-1)^{r}\overline{E}_{m}\biggl(-\frac{a(r+x)}{c}+y\biggl)\overline{B}_{n}\biggl(-\frac{b(r+x)}{c}+z\biggl)\nonumber\\
&&=\frac{nb^{n-1}}{2}\sum_{j=0}^{m}
\binom{m}{j}\frac{(-1)^{j}a^{m-j}}{c^{m+n-1-j}}\sum_{l=1}^{b}\overline{E}_{j}\biggl(\frac{a(l+z)}{b}-y\biggl)\nonumber\\
&&\qquad\times\overline{E}_{m+n-1-j}\biggl(\frac{c(l+z)}{b}-x\biggl)\nonumber\\
&&\quad+a^{m}\sum_{j=1}^{n}\binom{n}{j}\frac{(-1)^{j}b^{n-j}}{c^{m+n-j}}\sum_{l=1}^{a}(-1)^{l}
\overline{B}_{j}\biggl(\frac{b(l+y)}{a}-z\biggl)\nonumber\\
&&\qquad\times\overline{E}_{m+n-j}\biggl(\frac{c(l+y)}{a}-x\biggl)\nonumber\\
&&\quad+\frac{b^{n}}{a^{n}}\sum_{r=1}^{c}(-1)^{r}\overline{E}_{m+n}\biggl(-\frac{a(r+x)}{c}+y\biggl)\nonumber\\
&&\quad+\frac{1}{2}\delta_{0,m}\delta_{1,n}\sum_{r=1}^{c}(-1)^{-\frac{a(r+x)}{c}+y+r}\delta_{\mathbb{Z}}\biggl(-\frac{a(r+x)}{c}+y\biggl)
\delta_{\mathbb{Z}}\biggl(-\frac{b(r+x)}{c}+z\biggl),
\end{eqnarray*}
as was to be shown.
\end{proof}

In particular, we have the following result.

\begin{corollary}\label{cor2.7} Let $a,b,c\in\mathbb{N}$ with $2\nmid a$, $2\nmid c$ and $(a,b)=(b,c)=(a,c)=1$. Then
\begin{equation}\label{eq2.16}
s_{5}(a,b,c)+s_{5}(c,b,a)-\frac{1}{2}s_{4}(a,c,b)=\frac{1}{2}-\frac{b}{2ac},
\end{equation}
where $s_{4}(a,b,c)$, $s_{5}(a,b,c)$ are given for $a,b\in\mathbb{Z}$, $c\in\mathbb{N}$ by
\begin{equation*}
s_{4}(a,b,c)=\sum_{r=1}^{c-1}(-1)^{[\frac{ar}{c}]+[\frac{br}{c}]},
\end{equation*}
and
\begin{equation*}
s_{5}(a,b,c)=\sum_{r=1}^{c}(-1)^{r+[\frac{ar}{c}]}\biggl(\biggl(\frac{br}{c}\biggl)\biggl).
\end{equation*}
\end{corollary}

\begin{proof}
Taking $m=x=y=z=0$ and $n=1$ in Corollary \ref{cor2.6}, we see that for $2\nmid a$, $2\nmid c$ with $(a,b)=(b,c)=(a,c)=1$,
\begin{equation}\label{eq2.17}
s_{5}(a,b,c)+s_{5}(c,b,a)-\frac{1}{2}s_{4}(a,c,b)=\frac{1}{2}+\frac{b}{a}\sum_{r=1}^{c}(-1)^{r}\overline{E}_{1}\biggl(-\frac{ar}{c}\biggl).
\end{equation}
According to \eqref{eq2.2} and the familiar geometric sum, we find that for $n\in\mathbb{N}_{0}$, $a\in\mathbb{Z}\setminus\{0\}$, $c\in\mathbb{N}$, $x\in\mathbb{R}$ with $2\nmid a$ and $2\nmid c$,
\begin{eqnarray*}
\sum_{r=1}^{c}(-1)^{r}\overline{E}_{n}\biggl(x+\frac{ar}{c}\biggl)
&=&\frac{2n!}{(2\pi\mathrm{i})^{n+1}}\sum_{k=-\infty}^{+\infty}\frac{e^{2\pi\mathrm{i}(k-\frac{1}{2})x}}{\bigl(k-\frac{1}{2}\bigl)^{n+1}}
\sum_{r=1}^{c}e^{\frac{2\pi\mathrm{i}(ak-\frac{a-c}{2})r}{c}}\nonumber\\
&=&\frac{2cn!}{(2\pi\mathrm{i})^{n+1}}\sum_{\substack{k=-\infty\\c\mid(ak-\frac{a-c}{2})}}^{+\infty}\frac{e^{2\pi\mathrm{i}(k-\frac{1}{2})x}}{\bigl(k-\frac{1}{2}\bigl)^{n+1}}.
\end{eqnarray*}
For $k\in\mathbb{Z}$, the divisibility $c\mid\bigl(ak-\frac{a-c}{2}\bigr)$ implies that there exists $q\in\mathbb{Z}$ such that $a(2k-1)=c(2q-1)$. Hence, if $(a,c)=1$, then $c\mid(2k-1)$. It follows that for $n\in\mathbb{N}_{0}$, $a\in\mathbb{Z}\setminus\{0\}$, $c\in\mathbb{N}$, $x\in\mathbb{R}$ with $2\nmid a$, $2\nmid c$ and $(a,c)=1$,
\begin{eqnarray}\label{eq2.18}
\sum_{r=1}^{c}(-1)^{r}\overline{E}_{n}\biggl(x+\frac{ar}{c}\biggl)
&=&\frac{2n!}{(2\pi\mathrm{i})^{n+1}c^{n}}\sum_{k=-\infty}^{+\infty}\frac{e^{2\pi\mathrm{i}(k-\frac{1}{2})cx}}{\bigl(k-\frac{1}{2}\bigl)^{n+1}}\nonumber\\
&=&c^{-n}\overline{E}_{n}(cx).
\end{eqnarray}
Thus, applying \eqref{eq2.18} to the right-hand side of \eqref{eq2.17} gives the desired result.
\end{proof}

It is easily seen that the case $b=1$ and $c=1$ in Corollary \ref{cor2.7} gives the reciprocity formula \eqref{eq1.6} and \eqref{eq1.5}, respectively. In addition, swapping $a$
and $b$ in Corollary \ref{cor2.7}, and then using Pettet and Sitaramachandrarao's results stated in \cite[Equations (3.7), (3.8) and (3.9)]{pettet}, we rediscover Goldberg's \cite[Theorem 5.6]{goldberg} another three-term reciprocity formula
\begin{equation}\label{eq2.19}
s_{3}(ab',c)-\frac{1}{2}s_{4}(bc',a)+s_{1}(ca',b)=\frac{1}{2}-\frac{a}{2bc},
\end{equation}
where $a,b,c\in\mathbb{N}$ with $2\mid a$ and $(a,b)=(b,c)=(a,c)=1$, $a',b',c'\in\mathbb{Z}$ satisfy $2\mid a'$ and $aa'\equiv 1$ (mod $b$), $bb'\equiv 1$ (mod $c$), $cc'\equiv 1$ (mod $2a$).

\begin{theorem}\label{thm2.8} Let $m,n\in\mathbb{N}_{0}$, $a,b\in\mathbb{Z}\setminus\{0\}$ with $2\nmid a$ and $2\nmid b$. Then, for $x,y,z\in\mathbb{R}$,
\begin{eqnarray}\label{eq2.20}
&&\overline{E}_{m}(ax+y)\overline{E}_{n}(bx+z)\nonumber\\
&&=2b^{n}\mathrm{sgn}(b)
\sum_{j=0}^{m}\binom{m}{j}\frac{(-1)^{j}a^{m-j}}{m+n+1-j}\sum_{l=1}^{|b|}(-1)^{l}\overline{E}_{j}\biggl(\frac{a(l+z)}{b}-y\biggl)\nonumber\\
&&\qquad\times\overline{B}_{m+n+1-j}\biggl(x+\frac{l+z}{b}\biggl)\nonumber\\
&&\quad+2a^{m}\mathrm{sgn}(a)\sum_{j=0}^{n}\binom{n}{j}\frac{(-1)^{j}b^{n-j}}{m+n+1-j}\sum_{l=1}^{|a|}(-1)^{l}\overline{E}_{j}\biggl(\frac{b(l+y)}{a}-z\biggl)\nonumber\\
&&\qquad\times\overline{B}_{m+n+1-j}\biggl(x+\frac{l+y}{a}\biggl)\nonumber\\
&&\quad+\frac{(-1)^{n+1}2m!n!(a,b)^{m+n+2}}{a^{n+1}b^{m+1}(m+n+1)!}\overline{E}_{m+n+1}\biggl(\frac{by-az}{(a,b)}\biggl)\nonumber\\
&&\quad-\delta_{0,m}\delta_{0,n}\mathrm{sgn}(ab)(-1)^{ax+bx+y+z}\delta_{\mathbb{Z}}(ax+y)\delta_{\mathbb{Z}}(bx+z),
\end{eqnarray}
where $\delta_{\mathbb{Z}}(x)$ is as in \eqref{eq2.3}.
\end{theorem}

Clearly, taking $m=n=x=y=z=0$ in Theorem \ref{thm2.8}, we obtain
the reciprocity formula \eqref{eq1.6}.
We also another reciprocity formula among the generalized Hardy-Berndt sums \eqref{eq1.10} and \eqref{eq1.11}.

\begin{corollary}\label{cor2.9} Let $m,n\in\mathbb{N}_{0}$, $a,b,c\in\mathbb{N}$ with $2\nmid a$ and $2\nmid b$. Then, for $x,y,z\in\mathbb{R}$,
\begin{eqnarray}\label{eq2.21}
&&2b^{n}\sum_{j=0}^{m}\binom{m}{j}\frac{(-1)^{j}a^{m-j}}{c^{m+n-j}(m+n+1-j)}s_{j,m+n+1-j}^{(5)}\left(
\begin{matrix}
a & c & b \\
y & x & z
\end{matrix}
\right)\nonumber\\
&&\quad+2a^{m}\sum_{j=0}^{n}\binom{n}{j}\frac{(-1)^{j}b^{n-j}}{c^{m+n-j}(m+n+1-j)}s_{j,m+n+1-j}^{(5)}\left(
\begin{matrix}
b & c & a \\
z & x & y
\end{matrix}
\right)\nonumber\\
&&=s_{m,n}^{(4)}\left(
\begin{matrix}
-a & -b & c \\
-y & -z & x
\end{matrix}
\right)+\frac{(-1)^{n}2m!n!(a,b)^{m+n+2}c}{a^{n+1}b^{m+1}(m+n+1)!}\overline{E}_{m+n+1}\biggl(\frac{by-az}{(a,b)}\biggl)\nonumber\\
&&\quad+\delta_{0,m}\delta_{0,n}\sum_{r=1}^{c}(-1)^{\frac{(a+b)(r+x)}{c}-y-z}\nonumber\\
&&\qquad\times\delta_{\mathbb{Z}}\biggl(\frac{a(r+x)}{c}-y\biggl)
\delta_{\mathbb{Z}}\biggl(\frac{b(r+x)}{c}-z\biggl).
\end{eqnarray}
\end{corollary}

\begin{proof}
By replacing $x$ by $-(r+x)/c$ and then applying the operation $\sum_{r=1}^{c}$ to both sides of \eqref{eq2.20}, the desired result follows from \eqref{eq2.6}.
\end{proof}

It is trivial to see that taking $m=n=x=y=z=0$ and then swapping $b$ and $c$ in Corollary \ref{cor2.9}, we reobtain Corollary \ref{cor2.7}.

\begin{theorem}\label{thm2.10} Let $m,n\in\mathbb{N}_{0}$, $a,b\in\mathbb{Z}\setminus\{0\}$ with $2\nmid (a+b)$. Then, for $x,y,z\in\mathbb{R}$,
\begin{eqnarray}\label{eq2.22}
&&\overline{E}_{m}(ax+y)\overline{E}_{n}(bx+z)\nonumber\\
&&=b^{n}\mathrm{sgn}(b)
\sum_{j=0}^{m}\binom{m}{j}(-1)^{j+1}a^{m-j}\sum_{l=1}^{|b|}(-1)^{l}
\overline{E}_{j}\biggl(\frac{a(l+z)}{b}-y\biggl)\nonumber\\
&&\qquad\times\overline{E}_{m+n-j}\biggl(x+\frac{l+z}{b}\biggl)\nonumber\\
&&\quad+a^{m}\mathrm{sgn}(a)\sum_{j=0}^{n}\binom{n}{j}(-1)^{j+1}b^{n-j}
\sum_{l=1}^{|a|}(-1)^{l}\overline{E}_{j}\biggl(\frac{b(l+y)}{a}-z\biggl)\nonumber\\
&&\qquad\times\overline{E}_{m+n-j}\biggl(x+\frac{l+y}{a}\biggl)\nonumber\\
&&\quad-\delta_{0,m}\delta_{0,n}\mathrm{sgn}(ab)(-1)^{ax+bx+y+z}\delta_{\mathbb{Z}}(ax+y)\delta_{\mathbb{Z}}(bx+z),
\end{eqnarray}
where $\delta_{\mathbb{Z}}(x)$ is as in \eqref{eq2.3}.
\end{theorem}

If we take $m=n=x=y=z=0$ in Theorem \ref{thm2.10}, we obtain the reciprocity formula \eqref{eq1.3}.
As another application of Theorem \ref{thm2.10}, we have the following reciprocity formula for the generalized Hary-Berndt sums \eqref{eq1.7}.

\begin{corollary}\label{cor2.11} Let $m,n\in\mathbb{N}_{0}$, $a,b,c\in\mathbb{N}$ with $2\nmid (a+b)$ and $2\nmid c$. Then, for $x,y,z\in\mathbb{R}$,
\begin{eqnarray}\label{eq2.23}
&&b^{n}\sum_{j=0}^{m}\binom{m}{j}\frac{(-1)^{j+1}a^{m-j}}{c^{m+n-j}}S_{j,m+n-j}\left(
\begin{matrix}
a & c & b \\
y & x & z
\end{matrix}
\right)\nonumber\\
&&\quad+a^{m}\sum_{j=0}^{n}\binom{n}{j}\frac{(-1)^{j+1}b^{n-j}}{c^{m+n-j}}S_{j,m+n-j}\left(
\begin{matrix}
b & c & a \\
z & x & y
\end{matrix}
\right)\nonumber\\
&&=S_{m,n}\left(
\begin{matrix}
-a & -b & c \\
-y & -z & x
\end{matrix}
\right)\nonumber\\
&&\quad+\delta_{0,m}\delta_{0,n}\sum_{r=1}^{c}(-1)^{\frac{(a+b)(r+x)}{c}-y-z-r}\delta_{\mathbb{Z}}\biggl(\frac{a(r+x)}{c}-y\biggl)\nonumber\\
&&\qquad\times\delta_{\mathbb{Z}}\biggl(\frac{b(r+x)}{c}-z\biggl).
\end{eqnarray}
\end{corollary}

\begin{proof} Replacing $x$ by $-(r+x)/c$ and then applying the operation $\sum_{r=1}^{c}(-1)^{r}$ to both sides of \eqref{eq2.22}, we know from \eqref{eq2.15} that Corollary \ref{cor2.11} is true.
\end{proof}

It follows that we show the following extension of the reciprocity formula \eqref{eq1.3}.

\begin{corollary}\label{cor2.12} Let $a,b,c\in\mathbb{N}$ with $2\nmid (a+b)$, $2\nmid c$ and $(a,b)=(b,c)=(a,c)=1$. Then
\begin{equation}\label{eq2.24}
S(a,b,c)+S(b,c,a)+S(a,c,b)=1,
\end{equation}
where $S(a,b,c)$ is given for $a,b\in\mathbb{Z}$, $c\in\mathbb{N}$ by
\begin{equation*}
S(a,b,c)=\sum_{r=1}^{c-1}(-1)^{r+1+[\frac{ar}{c}]+[\frac{br}{c}]}.
\end{equation*}
\end{corollary}

\begin{proof}
Taking $m=n=x=y=z=0$ in Corollary \ref{cor2.11} gives the desired result.
\end{proof}

We remark that using Corollary \ref{cor2.12} and Pettet and Sitaramachandrarao's results stated in \cite[Equations (3.2), (3.3) and (3.4)]{pettet}, we rediscover Goldberg's \cite[Theorem 5.2]{goldberg} third three-term reciprocity formula
\begin{equation}\label{eq2.25}
S(bc',a)+S(ca',b)+S(ab',c)=1,
\end{equation}
where $a,b,c\in\mathbb{N}$ with $2\mid a$ and $(a,b)=(b,c)=(a,c)=1$, $a',b',c'\in\mathbb{Z}$ satisfy $2\mid a'$ and $aa'\equiv 1$ (mod $b$), $bb'\equiv 1$ (mod $c$), $cc'\equiv 1$ (mod $2a$).

\begin{remark}\label{rem2.13}
It is worth mentioning that Can \cite[Theorem 3]{can} in 2021 established a formula for products of two Bernoulli functions and used it to obtain some reciprocity formulas for the generalizations of the Hardy-Berndt sums involving the Bernoulli and Euler functions, where the $n$-th Bernoulli function is defined for $n\in\mathbb{N}_{0}$, $x\in\mathbb{R}$ by $\mathcal{B}_{n}(x)=B_{n}(\{x\})$ and the $n$-th Euler function is defined for $m\in\mathbb{Z}$, $n\in\mathbb{N}_{0}$, $0\leq x<1$ by $\mathcal{E}_{n}(x)=E_{n}(x)$ and $\mathcal{E}_{n}(x+m)=(-1)^{m}\mathcal{E}_{n}(x)$. For a further exploration, see the 2025 work of Isaacson \cite{isaacson}.

\end{remark}

\section{Some auxiliary results}

Before giving the proofs of Theorems \ref{thm2.1}, \ref{thm2.5}, \ref{thm2.8} and \ref{thm2.10}, we need the following auxiliary lemmas.

\begin{lemma}\label{lem3.1} (Parseval's formula) Suppose that $F(\theta)$ and $G(\theta)$ are two Riemann integrable, complex-valued functions on $\mathbb{R}$ of period $2\pi$ with the Fourier series
\begin{equation*}
F(\theta)=\sum_{n=-\infty}^{+\infty}a_{n}e^{\mathrm{i}n\theta},
\end{equation*}
and
\begin{equation*}
G(\theta)=\sum_{n=-\infty}^{+\infty}b_{n}e^{\mathrm{i}n\theta}.
\end{equation*}
Then
\begin{equation}\label{eq3.1}
\sum_{n=-\infty}^{+\infty}a_{n}\overline{b_{n}}=\frac{1}{2\pi}\int_{0}^{2\pi}F(\theta)\overline{G(\theta)}d\theta,
\end{equation}
where the horizontal bars indicate complex conjugation.
\end{lemma}

\begin{proof}
See \cite[Proposition 3.1.10]{cohen} or \cite[p. 81]{stein} for details.
\end{proof}

\begin{lemma}\label{lem3.2} Let $m,n\in\mathbb{N}$, $d,x,y\in\mathbb{R}$ be distinct. Then
\begin{eqnarray}\label{eq3.2}
\frac{1}{(d-x)^{m}(d-y)^{n}}&=&\sum_{j=1}^{m}\binom{m+n-j-1}{n-1}\frac{(-1)^{m-j}}{(x-y)^{m+n-j}(d-x)^{j}}\nonumber\\
&&+\sum_{j=1}^{n}\binom{m+n-j-1}{m-1}\frac{(-1)^{n-j}}{(y-x)^{m+n-j}(d-y)^{j}}.
\end{eqnarray}
\end{lemma}

\begin{proof}
See \cite[Lemma 2.2]{he2} for details.
\end{proof}

\begin{lemma}\label{lem3.3} Let $q,j\in\mathbb{N}$ with $q\geq2$, and let $\theta_{r}$ be a real-valued function defined on $r\in\mathbb{N}$ such that $\theta_{r}\not=0,\pm q,\pm 2q,\ldots$. Then
\begin{equation}\label{eq3.3}
\frac{\partial^{j-1}}{\partial a^{j-1}}\bigl(\cot(\pi a)\bigl)\biggl|_{a=\frac{\theta_{r}}{q}}=\frac{\delta_{1,j}}{\mathrm{i}}+2^{j}\pi^{j-1}\mathrm{i}^{j-2}F\biggl(\frac{\theta_{r}}{q},1-j\biggl),
\end{equation}
where $F(x,s)$ is the periodic-zeta function given for $x\in\mathbb{R}$, $s\in\mathbb{C}$ by
\begin{equation*}
F(x,s)=\sum_{n=1}^{\infty}\frac{e^{2\pi \mathrm{i}nx}}{n^{s}}\quad(\Re(s)>1).
\end{equation*}
\end{lemma}

\begin{proof}
See \cite[Equation (2.28)]{he1} for details.
\end{proof}

\begin{lemma}\label{lem3.4} Let $j\in\mathbb{N}$, $b,r\in\mathbb{Z}$ with $b\not=0$. Then
\begin{equation}\label{eq3.4}
\sum_{d=-\infty}^{+\infty}\frac{1}{\bigl(d+\frac{r}{b}-\frac{1}{2b}\bigl)^{j}}
=\frac{(-1)^{j}(2\pi\mathrm{i})^{j}\mathrm{sgn}(b)}{2(j-1)!b^{1-j}}\sum_{l=1}^{|b|}e^{\frac{2\pi \mathrm{i}lr}{b}-\frac{\pi \mathrm{i}l}{b}}\overline{E}_{j-1}\biggl(\frac{l}{b}\biggl).
\end{equation}
\end{lemma}

\begin{proof}
Since $\cot(a)$ has the following expression in partial fractions (see, e.g., \cite[p. 75]{abramowitz} or \cite[p. 327]{remmert})
\begin{eqnarray*}
\cot(a)&=&\frac{1}{a}+2a\sum_{d=1}^{\infty}\frac{1}{a^{2}-d^{2}\pi^{2}}\\
&=&\sum_{d=-\infty}^{+\infty}\frac{1}{a+d\pi}\quad(a\not=0,\pm \pi,\pm 2\pi,\ldots),
\end{eqnarray*}
we know that for $j\in\mathbb{N}$, $b,r\in\mathbb{Z}$ with $b\not=0$,
\begin{equation}\label{eq3.5}
\frac{\partial^{j-1}}{\partial a^{j-1}}\bigl(\cot(\pi a)\bigl)\biggl|_{a=\frac{r}{b}-\frac{1}{2b}}=\frac{(-1)^{j-1}(j-1)!}{\pi}\sum_{d=-\infty}^{+\infty}\frac{1}{\bigl(d+\frac{r}{b}-\frac{1}{2b}\bigl)^{j}}.
\end{equation}
By taking $\theta_{r}=q(\frac{r}{b}-\frac{1}{2b})$ in Lemma \ref{lem3.3}, we see from \eqref{eq3.5} that
\begin{equation}\label{eq3.6}
\sum_{d=-\infty}^{+\infty}\frac{1}{\bigl(d+\frac{r}{b}-\frac{1}{2b}\bigl)^{j}}=-\delta_{1,j}\pi\mathrm{i}+\frac{(-1)^{j}(2\pi\mathrm{i})^{j}}{(j-1)!}
F\biggl(\frac{r}{b}-\frac{1}{2b},1-j\biggl).
\end{equation}
It is clear that for $j\in\mathbb{N}_{0}$ (see, e.g., \cite[Equation (2.32)]{he3}),
\begin{equation}\label{eq3.7}
\eta(-j,x)=\frac{1}{2}E_{j}(x),
\end{equation}
where $\eta(s,x)$ is the alternating Hurwitz-zeta function given for $s\in\mathbb{C}$, $x\in\mathbb{R}$ with $x>0$ by
\begin{equation*}
\eta(s,x)=\sum_{n=0}^{\infty}(-1)^{n}\frac{1}{(n+x)^{s}}\quad(\Re(s)>0).
\end{equation*}
Hence, using the division algorithm and \eqref{eq3.7}, the periodic-zeta function on the right-hand side of \eqref{eq3.6} can be rewritten as
\begin{eqnarray}\label{eq3.8}
F\biggl(\frac{r}{b}-\frac{1}{2b},1-j\biggl)&=&\sum_{l=1}^{|b|}\sum_{m=0}^{\infty}\frac{e^{\frac{2\pi \mathrm{i}(m|b|+l)r}{b}}e^{-\frac{\pi \mathrm{i}(m|b|+l)}{b}}}{(m|b|+l)^{1-j}}\nonumber\\
&=&\sum_{l=1}^{|b|}e^{\frac{2\pi \mathrm{i}lr}{b}-\frac{\pi\mathrm{i}l}{b}}\sum_{m=0}^{\infty}(-1)^{m}\frac{1}{(m|b|+l)^{1-j}}\nonumber\\
&=&\frac{1}{2|b|^{1-j}}\sum_{l=1}^{|b|}e^{\frac{2\pi \mathrm{i}lr}{b}-\frac{\pi\mathrm{i}l}{b}}E_{j-1}\biggl(\frac{l}{|b|}\biggl).
\end{eqnarray}
Note that from $E_{0}(x)=1$ we obtain that for $j\in\mathbb{N}$, $l\in\mathbb{N}_{0}$, $b\in\mathbb{Z}$ with $1\leq l\leq |b|$,
\begin{equation}\label{eq3.9}
E_{j-1}\biggl(\frac{l}{|b|}\biggl)=\overline{E}_{j-1}\biggl(\frac{l}{|b|}\biggl)+\delta_{l,|b|}\delta_{1,j}.
\end{equation}
It follows from \eqref{eq3.8} and \eqref{eq3.9} that
\begin{equation}\label{eq3.10}
F\biggl(\frac{r}{b}-\frac{1}{2b},1-j\biggl)=\frac{1}{2|b|^{1-j}}\sum_{l=1}^{|b|}e^{\frac{2\pi \mathrm{i}lr}{b}-\frac{\pi\mathrm{i}l}{b}}\overline{E}_{j-1}\biggl(\frac{l}{|b|}\biggl)-\frac{1}{2}\delta_{1,j}.
\end{equation}
It is easily seen from \eqref{eq2.2} that for $n\in\mathbb{N}_{0}$, $x\in\mathbb{R}$,
\begin{equation}\label{eq3.11}
\overline{E}_{n}(-x)=\begin{cases}
\overline{E}_{n}(x),  &\text{if $x\in\mathbb{Z}$},\\
(-1)^{n+1}\overline{E}_{n}(x),  &\text{if $x\in\mathbb{R}\setminus\mathbb{Z}$}.
\end{cases}
\end{equation}
With the help of \eqref{eq1.13} and \eqref{eq3.11}, we rewrite \eqref{eq3.10} as
\begin{eqnarray}\label{eq3.12}
&&F\biggl(\frac{r}{b}-\frac{1}{2b},1-j\biggl)\nonumber\\
&&=\frac{\mathrm{sgn}(b)}{2b^{1-j}}\sum_{l=1}^{|b|}e^{\frac{2\pi \mathrm{i}lr}{b}-\frac{\pi\mathrm{i}l}{b}}\overline{E}_{j-1}\biggl(\frac{l}{b}\biggl)\nonumber\\
&&\quad+\frac{(\mathrm{sgn}(b))^{j-1}-\mathrm{sgn}(b)}{2b^{1-j}}\bigl(E_{j-1}(0)-\delta_{1,j}\bigl)-\frac{1}{2}\delta_{1,j}.
\end{eqnarray}
Since for $j\in\mathbb{N}$, (see, e.g., \cite[p. 805]{abramowitz})
\begin{equation}\label{eq3.13}
E_{j-1}(0)=-\frac{2(2^{j}-1)}{j}B_{j},
\end{equation}
and
\begin{equation}\label{eq3.14}
B_{2j+1}=0,
\end{equation}
where $B_{j}$ is the $j$-th Bernoulli number, by applying \eqref{eq3.13} and \eqref{eq3.14} to the right-hand side of \eqref{eq3.12}, we have
\begin{equation*}
F\biggl(\frac{r}{b}-\frac{1}{2b},1-j\biggl)=\frac{\mathrm{sgn}(b)}{2b^{1-j}}\sum_{l=1}^{|b|}e^{\frac{2\pi \mathrm{i}lr}{b}-\frac{\pi\mathrm{i}l}{b}}\overline{E}_{j-1}\biggl(\frac{l}{b}\biggl)-\frac{1}{2}\delta_{1,j},
\end{equation*}
from which and \eqref{eq3.6} we obtain \eqref{eq3.4} immediately.
This completes the proof of Lemma \ref{lem3.4}.
\end{proof}

\begin{lemma}\label{lem3.5}(Lerch's functional equation) If $0<x<1$ and $0<a\leq1$ then for $s\in\mathbb{C}$,
\begin{equation}\label{eq3.15}
\phi(x,a,1-s)=\frac{\Gamma(s)}{(2\pi)^{s}}\bigl(e^{\frac{\pi \mathrm{i}s}{2}-2\pi\mathrm{i}ax}\phi(-a,x,s)+e^{-\frac{\pi \mathrm{i}s}{2}+2\pi\mathrm{i}a(1-x)}\phi(a,1-x,s)\bigl),
\end{equation}
where $\phi(x,a,s)$ is the Lerch-zeta function given for $x,a\in\mathbb{R}$, $s\in\mathbb{C}$ with $0<a\leq1$ by
\begin{equation*}
\phi(x,a,s)=\sum_{n=0}^{\infty}\frac{e^{2\pi \mathrm{i}nx}}{(n+a)^{s}}\quad(\Re(s)>1).
\end{equation*}
\end{lemma}

\begin{proof}
See \cite{berndt1,lerch} for details.
\end{proof}

\begin{lemma}\label{lem3.6} Let $j\in\mathbb{N}$, $b,r\in\mathbb{Z}$ with $b\not=0$. Then, for $x\in\mathbb{R}$,
\begin{eqnarray}\label{eq3.16}
&&\sum_{d=-\infty}^{+\infty}\frac{e^{2\pi \mathrm{i}dx}}{\bigl(d+\frac{r}{b}-\frac{1}{2b}\bigl)^{j}}\nonumber\\
&&=\frac{(-1)^{j}(2\pi\mathrm{i})^{j}\mathrm{sgn}(b)}{2(j-1)!b^{1-j}}\sum_{l=1}^{|b|}e^{\frac{2\pi \mathrm{i}(l-x)r}{b}-\frac{\pi \mathrm{i}(l-x)}{b}}\overline{E}_{j-1}\biggl(\frac{l-x}{b}\biggl).
\end{eqnarray}
\end{lemma}

\begin{proof}
We first consider the case $x\in\mathbb{Z}$.
Since $l-x$ runs over a complete residue system modulo $|b|$ as
$l$ does. Hence, we know from \eqref{eq2.2} and Lemma \ref{lem3.4} that Lemma \ref{lem3.6} is true in the case when $x\in\mathbb{Z}$. We next discuss the case $x\not\in\mathbb{Z}$. It suffices to prove that Lemma \ref{lem3.6} holds in the case when $0<x<1$. Obviously, in this case, we have
\begin{eqnarray}\label{eq3.17}
&&\sum_{d=-\infty}^{+\infty}\frac{e^{2\pi \mathrm{i}dx}}{\bigl(d+\frac{r}{b}-\frac{1}{2b}\bigl)^{j}}\nonumber\\
&&=e^{-2\pi \mathrm{i}[\frac{r}{b}-\frac{1}{2b}]x}\biggl(\sum_{d=0}^{\infty}\frac{e^{2\pi \mathrm{i}dx}}{(d+\{\frac{r}{b}-\frac{1}{2b}\})^{j}}+(-1)^{j}e^{-2\pi \mathrm{i}x}\sum_{d=0}^{\infty}\frac{e^{-2\pi \mathrm{i}dx}}{(d+1-\{\frac{r}{b}-\frac{1}{2b}\})^{j}}\biggl)\nonumber\\
&&=e^{-2\pi \mathrm{i}[\frac{r}{b}-\frac{1}{2b}]x}\nonumber\\
&&\quad\times\biggl(\phi\biggl(x,\biggl\{\frac{r}{b}-\frac{1}{2b}\biggl\},j\biggl)
+(-1)^{j}e^{-2\pi \mathrm{i}x}\phi\biggl(-x,1-\biggl\{\frac{r}{b}-\frac{1}{2b}\biggl\},j\biggl)\biggl),
\end{eqnarray}
where $\phi(x,a,s)$ is as in \eqref{eq3.15}.
If we take $x=1-\{\frac{r}{b}-\frac{1}{2b}\}$ and $s=j$ and replace $a$ by $x$ in Lemma \ref{lem3.5}, then we have
\begin{eqnarray}\label{eq3.18}
&&\phi\biggl(1-\biggl\{\frac{r}{b}-\frac{1}{2b}\biggl\},x,1-j\biggl)\nonumber\\
&&=\frac{(j-1)!}{(2\pi)^{j}}\biggl(e^{\frac{\pi \mathrm{i}j}{2}-2\pi\mathrm{i}x(1-\{\frac{r}{b}-\frac{1}{2b}\})}\phi\biggl(-x,1-\biggl\{\frac{r}{b}-\frac{1}{2b}\biggl\},j\biggl)\nonumber\\
&&\quad+e^{-\frac{\pi \mathrm{i}j}{2}+2\pi\mathrm{i}x\{\frac{r}{b}-\frac{1}{2b}\}}\phi\biggl(x,\biggl\{\frac{r}{b}-\frac{1}{2b}\biggl\},j\biggl)\biggl).
\end{eqnarray}
Multiplying both sides of \eqref{eq3.18} by $e^{\frac{\pi \mathrm{i}j}{2}-2\pi\mathrm{i}\{\frac{r}{b}-\frac{1}{2b}\}x}$, and it follows from \eqref{eq3.17} that
\begin{equation}\label{eq3.19}
\sum_{d=-\infty}^{+\infty}\frac{e^{2\pi \mathrm{i}dx}}{\bigl(d+\frac{r}{b}-\frac{1}{2b}\bigl)^{j}}
=e^{-\frac{2\pi \mathrm{i}rx}{b}+\frac{\pi \mathrm{i}x}{b}}\frac{(2\pi \mathrm{i})^{j}\phi(1-\{\frac{r}{b}-\frac{1}{2b}\},x,1-j)}{(j-1)!}.
\end{equation}
Note that from the division algorithm, \eqref{eq1.13}, \eqref{eq3.7} and \eqref{eq3.11} we have
\begin{eqnarray}\label{eq3.20}
&&\phi\biggl(1-\biggl\{\frac{r}{b}-\frac{1}{2b}\biggl\},x,1-j\biggl)\nonumber\\
&&=\sum_{l=0}^{|b|-1}\sum_{m=0}^{\infty}\frac{e^{-\frac{2\pi \mathrm{i}(m|b|+l)r}{b}+\frac{\pi \mathrm{i}(m|b|+l)}{b}}}{(m|b|+l+x)^{1-j}}\nonumber\\
&&=\frac{1}{2|b|^{1-j}}\sum_{l=0}^{|b|-1}e^{-\frac{2\pi \mathrm{i}lr}{b}+\frac{\pi \mathrm{i}l}{b}}E_{j-1}\biggl(\frac{l+x}{|b|}\biggl)\nonumber\\
&&=\frac{\mathrm{sgn}(b)}{2b^{1-j}}\sum_{l=0}^{|b|-1}e^{-\frac{2\pi \mathrm{i}lr}{b}+\frac{\pi \mathrm{i}l}{b}}\overline{E}_{j-1}\biggl(\frac{l+x}{b}\biggl)\nonumber\\
&&=\frac{(-1)^{j}\mathrm{sgn}(b)}{2b^{1-j}}\sum_{l=1}^{|b|}e^{\frac{2\pi \mathrm{i}lr}{b}-\frac{\pi \mathrm{i}l}{b}}\overline{E}_{j-1}\biggl(\frac{l-x}{b}\biggl).
\end{eqnarray}
Therefore, inserting \eqref{eq3.20} into \eqref{eq3.19}, we say that Lemma \ref{lem3.6} holds in the case when $x\not\in\mathbb{Z}$. This completes the proof of Lemma \ref{lem3.6}.
\end{proof}

\begin{lemma}\label{lem3.7} Let $j\in\mathbb{N}$, $b,r\in\mathbb{Z}$ with $b\not=0$. Then, for $x\in\mathbb{R}$,
\begin{equation}\label{eq3.21}
\sideset{}{'}\sum_{d=-\infty}^{+\infty}\frac{e^{2\pi \mathrm{i}dx}}{\bigl(d+\frac{r}{b}\bigl)^{j}}
=\frac{(-1)^{j-1}(2\pi\mathrm{i})^{j}\mathrm{sgn}(b)}{j!b^{1-j}}\sum_{l=1}^{|b|}e^{\frac{2\pi \mathrm{i}(l-x)r}{b}}\overline{B}_{j}\biggl(\frac{l-x}{b}\biggl).
\end{equation}
\end{lemma}

\begin{proof}
See \cite[Lemma 2.7]{he2} for details.
\end{proof}

\begin{lemma}\label{lem3.8} Let $m\in\mathbb{N}_{0}$, $n\in\mathbb{N}$. Then
\begin{equation}\label{eq3.22}
\sum_{j=1}^{m+1}\binom{m+n-j}{n-1}=\binom{m+n}{n}.
\end{equation}
\end{lemma}

\begin{proof}
See \cite[Lemma 2.8]{he2} for details.
\end{proof}

\section{Proofs of Theorems}

\noindent{\em{The proof of Theorem \ref{thm2.1}.}}
Let
\begin{equation*}
f(x)=e^{\pi\mathrm{i}(ax+y)}\overline{E}_{m}(ax+y)\overline{B}_{n}(bx+z).
\end{equation*}
Since $f(x)$ is of bounded variation on every finite interval, $f(x)$ can be expanded in a Fourier series
\begin{equation}\label{eq4.1}
\frac{f(x^{+})+f(x^{-})}{2}=\sum_{k=-\infty}^{+\infty}c_{k}(m,n|a,b,y,z)e^{2\pi\mathrm{i}kx},
\end{equation}
where the Fourier coefficients $c_{k}(m,n|a,b,y,z)$ are given by
\begin{equation}\label{eq4.2}
c_{k}(m,n|a,b,y,z)=\int_{0}^{1}e^{\pi\mathrm{i}(ax+y)}\overline{E}_{m}(ax+y)\overline{B}_{n}(bx+z)e^{-2\pi\mathrm{i}kx}\text{d}x.
\end{equation}
If we make the change of the variable $\theta=2\pi x$ in Lemma \ref{lem3.1} and then set
\begin{equation*}
F(2\pi x)=e^{\pi\mathrm{i}(ax+y)}\overline{E}_{m}(ax+y),\quad\overline{G(2\pi x)}=\overline{B}_{n}(bx+z)e^{-2\pi\mathrm{i}kx},
\end{equation*}
then we find from \eqref{eq2.1}, \eqref{eq2.2}, \eqref{eq3.1} and \eqref{eq4.2} that
\begin{equation}\label{eq4.3}
c_{k}(m,n|a,b,y,z)=-\frac{2m!n!}{(2\pi \mathrm{i})^{m+n+1}}\underset{al+bj=k}{\sum_{l=-\infty}^{+\infty}\sideset{}{'}\sum_{j=-\infty}^{+\infty}}\frac{e^{2\pi\mathrm{i}(ly+jz)}} {(l-\frac{1}{2})^{m+1}j^{n}}.
\end{equation}
It is clear that the Diophantine equation $ax+by=k$ is solvable if and only if $(a,b)\mid k$, and in the case when it is solvable, all solutions of
$ax+by=k$ are given by
\begin{equation*}
x=\frac{k\overline{a}}{(a,b)}+\frac{bd}{(a,b)},\quad y=\frac{k\overline{b}}{(a,b)}-\frac{ad}{(a,b)},
\end{equation*}
where $\overline{a},\overline{b},d\in\mathbb{Z}$ satisfy $a\overline{a}+b\overline{b}=(a,b)$. Hence, we obtain from \eqref{eq4.3} that if $(a,b)\nmid k$ then
\begin{equation}\label{eq4.4}
c_{k}(m,n|a,b,y,z)=0,
\end{equation}
and if $(a,b)\mid k$ then
\begin{eqnarray}\label{eq4.5}
&&c_{k}(m,n|a,b,y,z)\nonumber\\
&&=-\frac{2m!n!(a,b)^{m+n+1}e^{\frac{2\pi\mathrm{i}k(\overline{a}y+\overline{b}z)}{(a,b)}}}{(2\pi \mathrm{i})^{m+n+1}}\sideset{}{'}\sum_{d=-\infty}^{+\infty}
\frac{e^{\frac{2\pi\mathrm{i}d(by-az)}{(a,b)}}}{\bigl(k\overline{a}+bd-\frac{(a,b)}{2}\bigl)^{m+1}
\bigl(k\overline{b}-ad\bigl)^{n}}\nonumber\\
&&=\frac{(-1)^{n+1}2m!n!(a,b)^{m+n+1}e^{\frac{2\pi\mathrm{i}k(\overline{a}y+\overline{b}z)}{(a,b)}}}{a^{n}b^{m+1}(2\pi \mathrm{i})^{m+n+1}}\nonumber\\
&&\qquad\times\sideset{}{'}\sum_{d=-\infty}^{+\infty}
\frac{e^{\frac{2\pi\mathrm{i}d(by-az)}{(a,b)}}}{\bigl(d+\frac{k\overline{a}}{b}-\frac{(a,b)}{2b}\bigl)^{m+1}
\bigl(d-\frac{k\overline{b}}{a}\bigl)^{n}}.
\end{eqnarray}
Since $2\nmid \overline{b}$ in the case when $2\mid a$ and $2\nmid b$, by taking $k=\frac{a}{2}$ in \eqref{eq4.5}, in light of \eqref{eq2.2}, we have
\begin{eqnarray}\label{eq4.6}
&&c_{\frac{a}{2}}(m,n|a,b,y,z)\nonumber\\
&&=\frac{(-1)^{n+1}2m!n!(a,b)^{m+n+1}e^{\frac{\pi\mathrm{i}a(\overline{a}y+\overline{b}z)}{(a,b)}}}{a^{n}b^{m+1}(2\pi\mathrm{i})^{m+n+1}}\sum_{d=-\infty}^{+\infty}
\frac{e^{\frac{2\pi\mathrm{i}d(by-az)}{(a,b)}}}{\bigl(d-\frac{\overline{b}-1}{2}-\frac{1}{2}\bigl)^{m+n+1}}\nonumber\\
&&=\frac{(-1)^{n+1}2m!n!(a,b)^{m+n+1}e^{\pi\mathrm{i}y}}{a^{n}b^{m+1}(2\pi\mathrm{i})^{m+n+1}}\sum_{d=-\infty}^{+\infty}
\frac{e^{\frac{2\pi\mathrm{i}(d-\frac{1}{2})(by-az)}{(a,b)}}}{\bigl(d-\frac{1}{2}\bigl)^{m+n+1}}\nonumber\\
&&=\frac{(-1)^{n+1}m!n!(a,b)^{m+n+1}e^{\pi\mathrm{i}y}}{a^{n}b^{m+1}(m+n)!}\overline{E}_{m+n}\biggl(\frac{by-az}{(a,b)}\biggl).
\end{eqnarray}
We next consider the case $k\not=\frac{a}{2}$ in \eqref{eq4.5}. Note that there is an implicit restriction
\begin{equation*}
\biggl(\frac{a}{(a,b)},\overline{b}\biggl)=\biggl(\frac{b}{(a,b)},\overline{a}\biggl)=1.
\end{equation*}
So if we take
\begin{equation*}
x=-\frac{k\overline{a}}{b}+\frac{(a,b)}{2b}, \quad y=\frac{k\overline{b}}{a}
\end{equation*}
in Lemma \ref{lem3.2}, then we see from \eqref{eq4.5} that
\begin{eqnarray}\label{eq4.7}
&&c_{k}(m,n|a,b,y,z)\nonumber\\
&&=\frac{(-1)^{n+1}2m!n!(a,b)^{m+n+1}e^{\frac{2\pi\mathrm{i}k(\overline{a}y+\overline{b}z)}{(a,b)}}}{a^{n}b^{m+1}(2\pi \mathrm{i})^{m+n+1}}\nonumber\\
&&\quad\times\biggl(\sum_{j=1}^{m+1}\binom{m+n-j}{n-1}(-1)^{m+1-j}\biggl(-\frac{ab}{k(a,b)-\frac{a}{2}(a,b)}\biggl)^{m+n+1-j}\nonumber\\
&&\qquad\quad\times\biggl(\sum_{d=-\infty}^{+\infty}
\frac{e^{\frac{2\pi\mathrm{i}d(by-az)}{(a,b)}}}{\bigl(d+\frac{k\overline{a}}{b}-\frac{(a,b)}{2b}\bigl)^{j}}-
\epsilon(a,k)e^{\frac{2\pi\mathrm{i}k\overline{b}(by-az)}{a(a,b)}}\biggl(\frac{ab}{k(a,b)-\frac{a}{2}(a,b)}\biggl)^{j}\biggl)\nonumber\\
&&\qquad+\sum_{j=1}^{n}\binom{m+n-j}{m}(-1)^{n-j}\biggl(\frac{ab}{k(a,b)-\frac{a}{2}(a,b)}\biggl)^{m+n+1-j}\nonumber\\
&&\qquad\quad\times\sideset{}{'}\sum_{d=-\infty}^{+\infty}\frac{e^{\frac{2\pi\mathrm{i}d(by-az)}{(a,b)}}}{\bigl(d-\frac{k\overline{b}}{a}\bigl)^{j}}\biggl),
\end{eqnarray}
where $\epsilon(a,k)=1$ or $0$ according to $a\mid k$ or $a\nmid k$. With the help of Lemmas \ref{lem3.6}, \ref{lem3.7} and \ref{lem3.8}, we discover that \eqref{eq4.7} can be rewritten as
\begin{eqnarray}\label{eq4.8}
&&c_{k}(m,n|a,b,y,z)\nonumber\\
&&=-b^{n-1}m!n!(a,b)\mathrm{sgn}(b)\sum_{j=1}^{m+1}\binom{m+n-j}{n-1}\frac{(-1)^{j}a^{m+1-j}}{(j-1)!\bigl(2\pi\mathrm{i}(k-\frac{a}{2})\bigl)^{m+n+1-j}}\nonumber\\
&&\qquad\times\sum_{l=1}^{|\frac{b}{(a,b)}|}e^{\frac{2\pi\mathrm{i}k(z+\overline{a}l)}{b}-\frac{\pi\mathrm{i}((a,b)l-by+az)}{b}}
\overline{E}_{j-1}\biggl(\frac{(a,b)l-by+az}{b}\biggl)\nonumber\\
&&\quad+2a^{m}m!n!(a,b)\mathrm{sgn}(a)\sum_{j=1}^{n}\binom{m+n-j}{m}\frac{(-1)^{j}b^{n-j}}{j!\bigl(2\pi\mathrm{i}(k-\frac{a}{2})\bigl)^{m+n+1-j}}\nonumber\\
&&\qquad\times\sum_{l=1}^{|\frac{a}{(a,b)}|}e^{\frac{2\pi\mathrm{i}k(y-\overline{b}l)}{a}}
\overline{B}_{j}\biggl(-\frac{(a,b)l-by+az}{a}\biggl)\nonumber\\
&&\quad+\epsilon(a,k)\frac{2a^{m+1}b^{n}(m+n)!e^{\frac{2\pi\mathrm{i}ky}{a}}}{\bigl(2\pi\mathrm{i}(k-\frac{a}{2})\bigl)^{m+n+1}}.
\end{eqnarray}
Inserting \eqref{eq4.4}, \eqref{eq4.6} and \eqref{eq4.8} into \eqref{eq4.1}, it follows from \eqref{eq2.2} that
\begin{eqnarray}\label{eq4.9}
&&e^{\pi\mathrm{i}(ax+y)}\overline{E}_{m}(ax+y)\overline{B}_{n}(bx+z)-\frac{1}{2}\delta_{0,m}\delta_{1,n}\mathrm{sgn}(ab)\delta_{\mathbb{Z}}(ax+y)\delta_{\mathbb{Z}}(bx+z)\nonumber\\
&&=\sideset{}{'}\sum_{k=-\infty}^{+\infty}\biggl(-b^{n-1}m!n!(a,b)\mathrm{sgn}(b)\sum_{j=1}^{m+1}
\binom{m+n-j}{n-1}\nonumber\\
&&\qquad\times\frac{(-1)^{j}a^{m+1-j}}{(j-1)!\bigl(2\pi\mathrm{i}((a,b)k-\frac{a}{2})\bigl)^{m+n+1-j}}\nonumber\\
&&\qquad\times\sum_{l=1}^{|\frac{b}{(a,b)}|}e^{\frac{2\pi\mathrm{i}(a,b)k(z+\overline{a}l)}{b}-\frac{\pi\mathrm{i}((a,b)l-by+az)}{b}}
\overline{E}_{j-1}\biggl(\frac{(a,b)l-by+az}{b}\biggl)\nonumber\\
&&\quad+2a^{m}m!n!(a,b)\mathrm{sgn}(a)\sum_{j=1}^{n}\binom{m+n-j}{m}\frac{(-1)^{j}b^{n-j}}{j!\bigl(2\pi\mathrm{i}((a,b)k-\frac{a}{2})\bigl)^{m+n+1-j}}\nonumber\\
&&\qquad\times\sum_{l=1}^{|\frac{a}{(a,b)}|}e^{\frac{2\pi\mathrm{i}(a,b)k(y-\overline{b}l)}{a}}
\overline{B}_{j}\biggl(-\frac{(a,b)l-by+az}{a}\biggl)\biggl)e^{2\pi\mathrm{i}(a,b)kx}\nonumber\\
&&\quad+\frac{(-1)^{n+1}m!n!(a,b)^{m+n+1}e^{\pi\mathrm{i}(ax+y)}}{a^{n}b^{m+1}(m+n)!}\overline{E}_{m+n}\biggl(\frac{by-az}{(a,b)}\biggl)\nonumber\\
&&\quad+\frac{b^{n}e^{\pi\mathrm{i}(ax+y)}}{a^{n}}\overline{E}_{m+n}(ax+y).
\end{eqnarray}
Since $\bigl(\frac{a}{(a,b)},\frac{b}{(a,b)}\bigl)=1$, $\frac{al}{(a,b)}$ runs over a complete residue
system modulo $|\frac{b}{(a,b)}|$ as $l$ does, and $-\frac{bl}{(a,b)}$ runs over a complete residue system modulo $|\frac{a}{(a,b)}|$
as $l$ does. It follows from \eqref{eq4.9} that
\begin{eqnarray*}
&&e^{\pi\mathrm{i}(ax+y)}\overline{E}_{m}(ax+y)\overline{B}_{n}(bx+z)-\frac{1}{2}\delta_{0,m}\delta_{1,n}\mathrm{sgn}(ab)\delta_{\mathbb{Z}}(ax+y)\delta_{\mathbb{Z}}(bx+z)\nonumber\\
&&=\sideset{}{'}\sum_{k=-\infty}^{+\infty}\biggl(-b^{n-1}m!n!(a,b)\mathrm{sgn}(b)\sum_{j=1}^{m+1}
\binom{m+n-j}{n-1}\nonumber\\
&&\qquad\times\frac{(-1)^{j}a^{m+1-j}}{(j-1)!\bigl(2\pi\mathrm{i}((a,b)k-\frac{a}{2})\bigl)^{m+n+1-j}}\nonumber\\
&&\qquad\times\sum_{l=1}^{|\frac{b}{(a,b)}|}e^{\frac{2\pi\mathrm{i}(a,b)k(l+z)}{b}-\frac{\pi\mathrm{i}(a(l+z)-by)}{b}}
\overline{E}_{j-1}\biggl(\frac{a(l+z)}{b}-y\biggl)\nonumber\\
&&\quad+2a^{m}m!n!(a,b)\mathrm{sgn}(a)\sum_{j=1}^{n}\binom{m+n-j}{m}\frac{(-1)^{j}b^{n-j}}{j!\bigl(2\pi\mathrm{i}((a,b)k-\frac{a}{2})\bigl)^{m+n+1-j}}\nonumber\\
&&\qquad\times\sum_{l=1}^{|\frac{a}{(a,b)}|}e^{\frac{2\pi\mathrm{i}(a,b)k(l+y)}{a}}
\overline{B}_{j}\biggl(\frac{b(l+y)}{a}-z\biggl)\biggl)e^{2\pi\mathrm{i}(a,b)kx}\nonumber\\
&&\quad+\frac{(-1)^{n+1}m!n!(a,b)^{m+n+1}e^{\pi\mathrm{i}(ax+y)}}{a^{n}b^{m+1}(m+n)!}\overline{E}_{m+n}\biggl(\frac{by-az}{(a,b)}\biggl)\nonumber\\
&&\quad+\frac{b^{n}e^{\pi\mathrm{i}(ax+y)}}{a^{n}}\overline{E}_{m+n}(ax+y),
\end{eqnarray*}
which together with $2(a,b)\mid a$ yields
\begin{eqnarray}\label{eq4.10}
&&\overline{E}_{m}(ax+y)\overline{B}_{n}(bx+z)-\frac{1}{2}\delta_{0,m}\delta_{1,n}\mathrm{sgn}(ab)(-1)^{ax+y}\delta_{\mathbb{Z}}(ax+y)\delta_{\mathbb{Z}}(bx+z)\nonumber\\
&&=\sideset{}{'}\sum_{k=-\infty}^{+\infty}\biggl(-b^{n-1}m!n!\mathrm{sgn}(b)\sum_{j=1}^{m+1}
\binom{m+n-j}{n-1}\nonumber\\
&&\qquad\times\frac{(-1)^{j}a^{m+1-j}}{(j-1)!(a,b)^{m+n-j}(2\pi\mathrm{i}k)^{m+n+1-j}}\nonumber\\
&&\qquad\times\sum_{l=1}^{|\frac{b}{(a,b)}|}e^{\frac{2\pi\mathrm{i}(a,b)k(l+z)}{b}}
\overline{E}_{j-1}\biggl(\frac{a(l+z)}{b}-y\biggl)\nonumber\\
&&\quad+2a^{m}m!n!\mathrm{sgn}(a)\sum_{j=1}^{n}\binom{m+n-j}{m}\frac{(-1)^{j}b^{n-j}}{j!(a,b)^{m+n-j}(2\pi\mathrm{i}k)^{m+n+1-j}}\nonumber\\
&&\qquad\times\sum_{l=1}^{|\frac{a}{(a,b)}|}(-1)^{l}e^{\frac{2\pi\mathrm{i}(a,b)k(l+y)}{a}}
\overline{B}_{j}\biggl(\frac{b(l+y)}{a}-z\biggl)\biggl)e^{2\pi\mathrm{i}(a,b)kx}\nonumber\\
&&\quad+\frac{(-1)^{n+1}m!n!(a,b)^{m+n+1}}{a^{n}b^{m+1}(m+n)!}\overline{E}_{m+n}\biggl(\frac{by-az}{(a,b)}\biggl)+\frac{b^{n}}{a^{n}}\overline{E}_{m+n}(ax+y).
\end{eqnarray}
Applying \eqref{eq2.1} to the right-hand side of \eqref{eq4.10}, we have
\begin{eqnarray}\label{eq4.11}
&&\overline{E}_{m}(ax+y)\overline{B}_{n}(bx+z)-\frac{1}{2}\delta_{0,m}\delta_{1,n}\mathrm{sgn}(ab)(-1)^{ax+y}\delta_{\mathbb{Z}}(ax+y)\delta_{\mathbb{Z}}(bx+z)\nonumber\\
&&=nb^{n-1}\mathrm{sgn}(b)
\sum_{j=1}^{m+1}\binom{m}{j-1}\frac{(-1)^{j}a^{m+1-j}}{(a,b)^{m+n-j}(m+n+1-j)}\nonumber\\
&&\qquad\times\sum_{l=1}^{|\frac{b}{(a,b)}|}
\overline{E}_{j-1}\biggl(\frac{a(l+z)}{b}-y\biggl)\overline{B}_{m+n+1-j}\biggl((a,b)x+\frac{(a,b)(l+z)}{b}\biggl)\nonumber\\
&&\quad-2a^{m}\mathrm{sgn}(a)\sum_{j=1}^{n}\binom{n}{j}\frac{(-1)^{j}b^{n-j}}{(a,b)^{m+n-j}(m+n+1-j)}\nonumber\\
&&\qquad\times\sum_{l=1}^{|\frac{a}{(a,b)}|}(-1)^{l}\overline{B}_{j}\biggl(\frac{b(l+y)}{a}-z\biggl)\overline{B}_{m+n+1-j}\biggl((a,b)x+\frac{(a,b)(l+y)}{a}\biggl)\nonumber\\
&&\quad+\frac{(-1)^{n+1}m!n!(a,b)^{m+n+1}}{a^{n}b^{m+1}(m+n)!}\overline{E}_{m+n}\biggl(\frac{by-az}{(a,b)}\biggl)+\frac{b^{n}}{a^{n}}\overline{E}_{m+n}(ax+y).
\end{eqnarray}
Since $2\mid a$ and $2\nmid b$, by \eqref{eq1.12}, \eqref{eq1.13}, \eqref{eq2.5}, the property of residue systems and the division algorithm, we have
\begin{eqnarray}\label{eq4.12}
&&\frac{1}{(a,b)^{m+n-j}}\sum_{l=1}^{|\frac{b}{(a,b)}|}
\overline{E}_{j-1}\biggl(\frac{a(l+z)}{b}-y\biggl)\nonumber\\
&&\qquad\times\overline{B}_{m+n+1-j}\biggl((a,b)x+\frac{(a,b)(l+z)}{b}\biggl)\nonumber\\
&&=\sum_{k=0}^{(a,b)-1}\sum_{l=1}^{|\frac{b}{(a,b)}|}\overline{E}_{j-1}\biggl(\frac{a\bigl(k|\frac{b}{(a,b)}|+l+z\bigl)}{b}-y\biggl)\nonumber\\
&&\qquad\times\overline{B}_{m+n+1-j}\biggl(x+\frac{l+z+k|\frac{b}{(a,b)}|}{b}\biggl)\nonumber\\
&&=\sum_{l=1}^{|b|}\overline{E}_{j-1}\biggl(\frac{a(l+z)}{b}-y\biggl)\overline{B}_{m+n+1-j}\biggl(x+\frac{l+z}{b}\biggl),
\end{eqnarray}
and
\begin{eqnarray}\label{eq4.13}
&&\frac{1}{(a,b)^{m+n-j}}\sum_{l=1}^{|\frac{a}{(a,b)}|}(-1)^{l}\overline{B}_{j}\biggl(\frac{b(l+y)}{a}-z\biggl)\nonumber\\
&&\qquad\times\overline{B}_{m+n+1-j}\biggl((a,b)x+\frac{(a,b)(l+y)}{a}\biggl)\nonumber\\
&&=\sum_{k=0}^{(a,b)-1}\sum_{l=1}^{|\frac{a}{(a,b)}|}(-1)^{k|\frac{a}{(a,b)}|+l}\overline{B}_{j}\biggl(\frac{b\bigl(k|\frac{a}{(a,b)}|+l+y\bigl)}{a}-z\biggl)\nonumber\\
&&\qquad\times\overline{B}_{m+n+1-j}\biggl(x+\frac{l+y+k|\frac{a}{(a,b)}|}{a}\biggl)\nonumber\\
&&=\sum_{l=1}^{|a|}(-1)^{l}\overline{B}_{j}\biggl(\frac{b(l+y)}{a}-z\biggl)\overline{B}_{m+n+1-j}\biggl(x+\frac{l+y}{a}\biggl).
\end{eqnarray}
Thus, by inserting \eqref{eq4.12} and \eqref{eq4.13} into \eqref{eq4.11}, we get \eqref{eq2.3} and finish the proof of Theorem \ref{thm2.1}.

\noindent{\em{The proof of Theorem \ref{thm2.5}.}}
Since $2\nmid a$, by inserting \eqref{eq4.4} and \eqref{eq4.8} into \eqref{eq4.1}, we obtain that
\begin{eqnarray*}
&&e^{\pi\mathrm{i}(ax+y)}\overline{E}_{m}(ax+y)\overline{B}_{n}(bx+z)-\frac{1}{2}\delta_{0,m}\delta_{1,n}\mathrm{sgn}(ab)\delta_{\mathbb{Z}}(ax+y)\delta_{\mathbb{Z}}(bx+z)\nonumber\\
&&=\sum_{k=-\infty}^{+\infty}\biggl(-b^{n-1}m!n!(a,b)\mathrm{sgn}(b)\sum_{j=1}^{m+1}
\binom{m+n-j}{n-1}\nonumber\\
&&\qquad\times\frac{(-1)^{j}a^{m+1-j}}{(j-1)!\bigl(2\pi\mathrm{i}((a,b)k-\frac{a}{2})\bigl)^{m+n+1-j}}\nonumber\\
&&\qquad\times\sum_{l=1}^{|\frac{b}{(a,b)}|}e^{\frac{2\pi\mathrm{i}(a,b)k(z+\overline{a}l)}{b}-\frac{\pi\mathrm{i}((a,b)l-by+az)}{b}}
\overline{E}_{j-1}\biggl(\frac{(a,b)l-by+az}{b}\biggl)\nonumber\\
&&\quad+2a^{m}m!n!(a,b)\mathrm{sgn}(a)\sum_{j=1}^{n}\binom{m+n-j}{m}\frac{(-1)^{j}b^{n-j}}{j!\bigl(2\pi\mathrm{i}((a,b)k-\frac{a}{2})\bigl)^{m+n+1-j}}\nonumber\\
&&\qquad\times\sum_{l=1}^{|\frac{a}{(a,b)}|}e^{\frac{2\pi\mathrm{i}(a,b)k(y-\overline{b}l)}{a}}
\overline{B}_{j}\biggl(-\frac{(a,b)l-by+az}{a}\biggl)\biggl)e^{2\pi\mathrm{i}(a,b)kx}\nonumber\\
&&\quad+\frac{b^{n}e^{\pi\mathrm{i}(ax+y)}}{a^{n}}\overline{E}_{m+n}(ax+y),
\end{eqnarray*}
from which and the properties of residue system we have
\begin{eqnarray*}
&&e^{\pi\mathrm{i}(ax+y)}\overline{E}_{m}(ax+y)\overline{B}_{n}(bx+z)-\frac{1}{2}\delta_{0,m}\delta_{1,n}\mathrm{sgn}(ab)\delta_{\mathbb{Z}}(ax+y)\delta_{\mathbb{Z}}(bx+z)\nonumber\\
&&=\sum_{k=-\infty}^{+\infty}\biggl(-b^{n-1}m!n!(a,b)\mathrm{sgn}(b)\sum_{j=1}^{m+1}
\binom{m+n-j}{n-1}\nonumber\\
&&\qquad\times\frac{(-1)^{j}a^{m+1-j}}{(j-1)!\bigl(2\pi\mathrm{i}((a,b)k-\frac{a}{2})\bigl)^{m+n+1-j}}\nonumber\\
&&\qquad\times\sum_{l=1}^{|\frac{b}{(a,b)}|}e^{\frac{2\pi\mathrm{i}(a,b)k(l+z)}{b}-\frac{\pi\mathrm{i}(a(l+z)-by)}{b}}
\overline{E}_{j-1}\biggl(\frac{a(l+z)}{b}-y\biggl)\nonumber\\
&&\quad+2a^{m}m!n!(a,b)\mathrm{sgn}(a)\sum_{j=1}^{n}\binom{m+n-j}{m}\frac{(-1)^{j}b^{n-j}}{j!\bigl(2\pi\mathrm{i}((a,b)k-\frac{a}{2})\bigl)^{m+n+1-j}}\nonumber\\
&&\qquad\times\sum_{l=1}^{|\frac{a}{(a,b)}|}e^{\frac{2\pi\mathrm{i}(a,b)k(l+y)}{a}}
\overline{B}_{j}\biggl(\frac{b(l+y)}{a}-z\biggl)\biggl)e^{2\pi\mathrm{i}(a,b)kx}\nonumber\\
&&\quad+\frac{b^{n}e^{\pi\mathrm{i}(ax+y)}}{a^{n}}\overline{E}_{m+n}(ax+y).
\end{eqnarray*}
It follows that
\begin{eqnarray}\label{eq4.14}
&&\overline{E}_{m}(ax+y)\overline{B}_{n}(bx+z)-\frac{1}{2}\delta_{0,m}\delta_{1,n}\mathrm{sgn}(ab)(-1)^{ax+y}\delta_{\mathbb{Z}}(ax+y)\delta_{\mathbb{Z}}(bx+z)\nonumber\\
&&=\sum_{k=-\infty}^{+\infty}\biggl(-b^{n-1}m!n!\mathrm{sgn}(b)\sum_{j=1}^{m+1}
\binom{m+n-j}{n-1}\nonumber\\
&&\qquad\times\frac{(-1)^{j}a^{m+1-j}}{(j-1)!(a,b)^{m+n-j}\bigl(2\pi\mathrm{i}(k-\frac{1}{2})\bigl)^{m+n+1-j}}\nonumber\\
&&\qquad\times\sum_{l=1}^{|\frac{b}{(a,b)}|}e^{\frac{2\pi\mathrm{i}(a,b)(k-\frac{1}{2})(l+z)}{b}}
\overline{E}_{j-1}\biggl(\frac{a(l+z)}{b}-y\biggl)\nonumber\\
&&\quad+2a^{m}m!n!\mathrm{sgn}(a)\sum_{j=1}^{n}\binom{m+n-j}{m}\nonumber\\
&&\qquad\times\frac{(-1)^{j}b^{n-j}}{j!(a,b)^{m+n-j}\bigl(2\pi\mathrm{i}(k-\frac{1}{2})\bigl)^{m+n+1-j}}\nonumber\\
&&\qquad\times\sum_{l=1}^{|\frac{a}{(a,b)}|}(-1)^{l}e^{\frac{2\pi\mathrm{i}(a,b)(k-\frac{1}{2})(l+y)}{a}}
\overline{B}_{j}\biggl(\frac{b(l+y)}{a}-z\biggl)\biggl)e^{2\pi\mathrm{i}(a,b)(k-\frac{1}{2})x}\nonumber\\
&&\quad+\frac{b^{n}}{a^{n}}\overline{E}_{m+n}(ax+y).
\end{eqnarray}
By applying \eqref{eq2.2} to the right-hand side of \eqref{eq4.14}, we have
\begin{eqnarray}\label{eq4.15}
&&\overline{E}_{m}(ax+y)\overline{B}_{n}(bx+z)-\frac{1}{2}\delta_{0,m}\delta_{1,n}\mathrm{sgn}(ab)(-1)^{ax+y}\delta_{\mathbb{Z}}(ax+y)\delta_{\mathbb{Z}}(bx+z)\nonumber\\
&&=-\frac{nb^{n-1}\mathrm{sgn}(b)}{2}\sum_{j=1}^{m+1}
\binom{m}{j-1}\frac{(-1)^{j}a^{m+1-j}}{(a,b)^{m+n-j}}\sum_{l=1}^{|\frac{b}{(a,b)}|}\overline{E}_{j-1}\biggl(\frac{a(l+z)}{b}-y\biggl)\nonumber\\
&&\qquad\times\overline{E}_{m+n-j}\biggl((a,b)x+\frac{(a,b)(l+z)}{b}\biggl)\nonumber\\
&&\quad+a^{m}\mathrm{sgn}(a)\sum_{j=1}^{n}\binom{n}{j}\frac{(-1)^{j}b^{n-j}}{(a,b)^{m+n-j}}\sum_{l=1}^{|\frac{a}{(a,b)}|}(-1)^{l}
\overline{B}_{j}\biggl(\frac{b(l+y)}{a}-z\biggl)\nonumber\\
&&\qquad\times\overline{E}_{m+n-j}\biggl((a,b)x+\frac{(a,b)(l+y)}{a}\biggl)\nonumber\\
&&\quad+\frac{b^{n}}{a^{n}}\overline{E}_{m+n}(ax+y).
\end{eqnarray}
Using the same method to prove \eqref{eq2.18}, we find that for $n\in\mathbb{N}_{0}$, $a\in\mathbb{Z}\setminus\{0\}$, $c\in\mathbb{N}$, $x\in\mathbb{R}$ with $2\nmid a$, $2\nmid c$ and $(a,c)=1$,
\begin{equation}\label{eq4.16}
\sum_{r=0}^{c-1}(-1)^{r}\overline{E}_{n}\biggl(x+\frac{ar}{c}\biggl)
=c^{-n}\overline{E}_{n}(cx).
\end{equation}
Since $2\nmid a$, by \eqref{eq1.12}, \eqref{eq1.13}, \eqref{eq4.16} and the division algorithm that
\begin{eqnarray}\label{eq4.17}
&&\frac{1}{(a,b)^{m+n-j}}\sum_{l=1}^{|\frac{b}{(a,b)}|}
\overline{E}_{j-1}\biggl(\frac{a(l+z)}{b}-y\biggl)\nonumber\\
&&\qquad\times\overline{E}_{m+n-j}\biggl((a,b)x+\frac{(a,b)(l+z)}{b}\biggl)\nonumber\\
&&=\sum_{k=0}^{(a,b)-1}\sum_{l=1}^{|\frac{b}{(a,b)}|}\overline{E}_{j-1}\biggl(\frac{a\bigl(k|\frac{b}{(a,b)}|+l+z\bigl)}{b}-y\biggl)\nonumber\\
&&\qquad\times\overline{E}_{m+n-j}\biggl(x+\frac{l+z+k|\frac{b}{(a,b)}|}{b}\biggl)\nonumber\\
&&=\sum_{l=1}^{|b|}\overline{E}_{j-1}\biggl(\frac{a(l+z)}{b}-y\biggl)\overline{E}_{m+n-j}\biggl(x+\frac{l+z}{b}\biggl),
\end{eqnarray}
and
\begin{eqnarray}\label{eq4.18}
&&\frac{1}{(a,b)^{m+n-j}}\sum_{l=1}^{|\frac{a}{(a,b)}|}(-1)^{l}\overline{B}_{j}\biggl(\frac{b(l+y)}{a}-z\biggl)\nonumber\\
&&\qquad\times\overline{E}_{m+n-j}\biggl((a,b)x+\frac{(a,b)(l+y)}{a}\biggl)\nonumber\\
&&=\sum_{k=0}^{(a,b)-1}\sum_{l=1}^{|\frac{a}{(a,b)}|}(-1)^{k|\frac{a}{(a,b)}|+l}\overline{B}_{j}\biggl(\frac{b\bigl(k|\frac{a}{(a,b)}|+l+y\bigl)}{a}-z\biggl)\nonumber\\
&&\qquad\times\overline{E}_{m+n-j}\biggl(x+\frac{l+y+k|\frac{a}{(a,b)}|}{a}\biggl)\nonumber\\
&&=\sum_{l=1}^{|a|}(-1)^{l}\overline{B}_{j}\biggl(\frac{b(l+y)}{a}-z\biggl)\overline{E}_{m+n-j}\biggl(x+\frac{l+y}{a}\biggl).
\end{eqnarray}
Therefore, by inserting \eqref{eq4.17} and \eqref{eq4.18} into \eqref{eq4.15}, we obtain \eqref{eq2.13} and prove Theorem \ref{thm2.5}.

\noindent{\em{The proof of Theorem \ref{thm2.8}.}}
Let
\begin{equation*}
g(x)=e^{\pi\mathrm{i}(ax+bx+y+z)}\overline{E}_{m}(ax+y)\overline{E}_{n}(bx+z).
\end{equation*}
Since $g(x)$ is of bounded variation on every finite interval, $g(x)$ may be expanded in a Fourier series
\begin{equation}\label{eq4.19}
\frac{g(x^{+})+g(x^{-})}{2}=\sum_{k=-\infty}^{+\infty}\widetilde{c}_{k}(m,n|a,b,y,z)e^{2\pi\mathrm{i}kx},
\end{equation}
where the Fourier coefficients $\widetilde{c}_{k}(m,n|a,b,y,z)$ are determined by
\begin{equation}\label{eq4.20}
\widetilde{c}_{k}(m,n|a,b,y,z)=\int_{0}^{1}e^{\pi\mathrm{i}(ax+y)}\overline{E}_{m}(ax+y)e^{\pi\mathrm{i}(bx+z)}\overline{E}_{n}(bx+z)e^{-2\pi\mathrm{i}kx}\text{d}x.
\end{equation}
By making the change of the variable $\theta=2\pi x$ in Lemma \ref{lem3.1} and setting
\begin{equation*}
F(2\pi x)=e^{\pi\mathrm{i}(ax+y)}\overline{E}_{m}(ax+y),\quad\overline{G(2\pi x)}=e^{\pi\mathrm{i}(bx+z)}\overline{E}_{n}(bx+z)e^{-2\pi\mathrm{i}kx},
\end{equation*}
we obtain from \eqref{eq2.2}, \eqref{eq3.1} and \eqref{eq4.20} that
\begin{equation}\label{eq4.21}
\widetilde{c}_{k}(m,n|a,b,y,z)=\frac{4m!n!}{(2\pi \mathrm{i})^{m+n+2}}\underset{al+bj=k}{\sum_{l=-\infty}^{+\infty}\sum_{j=-\infty}^{+\infty}}\frac{e^{2\pi\mathrm{i}(ly+jz)}} {(l-\frac{1}{2})^{m+1}(j-\frac{1}{2})^{n+1}}.
\end{equation}
In a similar consideration to \eqref{eq4.4} and \eqref{eq4.5}, we know from \eqref{eq4.21} that if $(a,b)\nmid k$ then
\begin{equation}\label{eq4.22}
\widetilde{c}_{k}(m,n|a,b,y,z)=0,
\end{equation}
and if $(a,b)\mid k$ then
\begin{eqnarray}\label{eq4.23}
&&\widetilde{c}_{k}(m,n|a,b,y,z)\nonumber\\
&&=\frac{4m!n!(a,b)^{m+n+2}e^{\frac{2\pi\mathrm{i}k(\overline{a}y+\overline{b}z)}{(a,b)}}}{(2\pi \mathrm{i})^{m+n+2}}\nonumber\\
&&\qquad\times\sum_{d=-\infty}^{+\infty}
\frac{e^{\frac{2\pi\mathrm{i}d(by-az)}{(a,b)}}}{\bigl(k\overline{a}+bd-\frac{(a,b)}{2}\bigl)^{m+1}
\bigl(k\overline{b}-ad-\frac{(a,b)}{2}\bigl)^{n+1}}\nonumber\\
&&=\frac{(-1)^{n+1}4m!n!(a,b)^{m+n+2}e^{\frac{2\pi\mathrm{i}k(\overline{a}y+\overline{b}z)}{(a,b)}}}{a^{n+1}b^{m+1}(2\pi \mathrm{i})^{m+n+2}}\nonumber\\
&&\qquad\times\sum_{d=-\infty}^{+\infty}
\frac{e^{\frac{2\pi\mathrm{i}d(by-az)}{(a,b)}}}{\bigl(d+\frac{k\overline{a}}{b}-\frac{(a,b)}{2b}\bigl)^{m+1}
\bigl(d-\frac{k\overline{b}}{a}+\frac{(a,b)}{2a}\bigl)^{n+1}}.
\end{eqnarray}
Since $\overline{a}$ and $\overline{b}$ have the opposite parity in the case when $2\nmid a$ and $2\nmid b$, by setting $k=\frac{a+b}{2}$ in \eqref{eq4.23}, in view of \eqref{eq2.2}, we have
\begin{eqnarray}\label{eq4.24}
&&\widetilde{c}_{\frac{a+b}{2}}(m,n|a,b,y,z)\nonumber\\
&&=\frac{(-1)^{n+1}4m!n!(a,b)^{m+n+2}e^{\frac{\pi\mathrm{i}(a+b)(\overline{a}y+\overline{b}z)}{(a,b)}}}{a^{n+1}b^{m+1}(2\pi\mathrm{i})^{m+n+2}}\sum_{d=-\infty}^{+\infty}
\frac{e^{\frac{2\pi\mathrm{i}d(by-az)}{(a,b)}}}{\bigl(d+\frac{\overline{a}-\overline{b}}{2}\bigl)^{m+n+2}}\nonumber\\
&&=\frac{(-1)^{n+1}4m!n!(a,b)^{m+n+2}e^{\frac{\pi\mathrm{i}(a+b)(\overline{a}y+\overline{b}z)}{(a,b)}}}{a^{n+1}b^{m+1}(2\pi\mathrm{i})^{m+n+2}}\sum_{d=-\infty}^{+\infty}
\frac{e^{\frac{2\pi\mathrm{i}d(by-az)}{(a,b)}}}{\bigl(d+\frac{\overline{a}-\overline{b}+1}{2}-\frac{1}{2}\bigl)^{m+n+2}}\nonumber\\
&&=\frac{(-1)^{n+1}4m!n!(a,b)^{m+n+2}e^{\pi\mathrm{i}(y+z)}}{a^{n+1}b^{m+1}(2\pi\mathrm{i})^{m+n+2}}\sum_{d=-\infty}^{+\infty}
\frac{e^{\frac{2\pi\mathrm{i}(d-\frac{1}{2})(by-az)}{(a,b)}}}{\bigl(d-\frac{1}{2}\bigl)^{m+n+2}}\nonumber\\
&&=\frac{(-1)^{n+1}2m!n!(a,b)^{m+n+2}e^{\pi\mathrm{i}(y+z)}}{a^{n+1}b^{m+1}(m+n+1)!}\overline{E}_{m+n+1}\biggl(\frac{by-az}{(a,b)}\biggl).
\end{eqnarray}
We next consider the case $k\not=\frac{a+b}{2}$ in \eqref{eq4.23}. By taking
\begin{equation*}
x=-\frac{k\overline{a}}{b}+\frac{(a,b)}{2b}, \quad y=\frac{k\overline{b}}{a}-\frac{(a,b)}{2a}
\end{equation*}
in Lemma \ref{lem3.2}, we know from \eqref{eq4.23} that
\begin{eqnarray}\label{eq4.25}
&&\widetilde{c}_{k}(m,n|a,b,y,z)\nonumber\\
&&=\frac{(-1)^{n+1}4m!n!(a,b)^{m+n+2}e^{\frac{2\pi\mathrm{i}k(\overline{a}y+\overline{b}z)}{(a,b)}}}{a^{n+1}b^{m+1}(2\pi \mathrm{i})^{m+n+2}}\nonumber\\
&&\quad\times\biggl(\sum_{j=1}^{m+1}\binom{m+n+1-j}{n}(-1)^{m+1-j}\biggl(-\frac{ab}{k(a,b)-(a,b)\frac{a+b}{2}}\biggl)^{m+n+2-j}\nonumber\\
&&\qquad\quad\times\sum_{d=-\infty}^{+\infty}
\frac{e^{\frac{2\pi\mathrm{i}d(by-az)}{(a,b)}}}{\bigl(d+\frac{k\overline{a}}{b}-\frac{(a,b)}{2b}\bigl)^{j}}\nonumber\\
&&\qquad+\sum_{j=1}^{n+1}\binom{m+n+1-j}{m}(-1)^{n+1-j}\biggl(\frac{ab}{k(a,b)-(a,b)\frac{a+b}{2}}\biggl)^{m+n+2-j}\nonumber\\
&&\qquad\quad\times\sum_{d=-\infty}^{+\infty}\frac{e^{\frac{2\pi\mathrm{i}d(by-az)}{(a,b)}}}{\bigl(d-\frac{k\overline{b}}{a}+\frac{(a,b)}{2a}\bigl)^{j}}\biggl).
\end{eqnarray}
Applying Lemma \ref{lem3.6} to the right-hand side of \eqref{eq4.25} and then using the properties of residue systems, we obtain from $(a,b)\mid k$ that
\begin{eqnarray}\label{eq4.26}
&&\widetilde{c}_{k}(m,n|a,b,y,z)\nonumber\\
&&=2b^{n}m!n!(a,b)\mathrm{sgn}(b)\sum_{j=1}^{m+1}\binom{m+n+1-j}{n}\frac{(-1)^{j}a^{m+1-j}}{(j-1)!\bigl(2\pi\mathrm{i}(k-\frac{a+b}{2})\bigl)^{m+n+2-j}}\nonumber\\
&&\qquad\times\sum_{l=1}^{|\frac{b}{(a,b)}|}e^{\frac{2\pi\mathrm{i}k(z+\overline{a}l)}{b}-\frac{\pi\mathrm{i}((a,b)l-by+az)}{b}}
\overline{E}_{j-1}\biggl(\frac{(a,b)l-by+az}{b}\biggl)\nonumber\\
&&\quad+2a^{m}m!n!(a,b)\mathrm{sgn}(a)\sum_{j=1}^{n+1}\binom{m+n+1-j}{m}\frac{(-1)^{j}b^{n+1-j}}{(j-1)!\bigl(2\pi\mathrm{i}(k-\frac{a+b}{2})\bigl)^{m+n+2-j}}\nonumber\\
&&\qquad\times\sum_{l=1}^{|\frac{a}{(a,b)}|}e^{\frac{2\pi\mathrm{i}k(y-\overline{b}l)}{a}+\frac{\pi\mathrm{i}((a,b)l-by+az)}{a}}
\overline{E}_{j-1}\biggl(-\frac{(a,b)l-by+az}{a}\biggl)\nonumber\\
&&=2b^{n}m!n!(a,b)\mathrm{sgn}(b)\sum_{j=1}^{m+1}\binom{m+n+1-j}{n}\frac{(-1)^{j}a^{m+1-j}}{(j-1)!\bigl(2\pi\mathrm{i}(k-\frac{a+b}{2})\bigl)^{m+n+2-j}}\nonumber\\
&&\qquad\times\sum_{l=1}^{|\frac{b}{(a,b)}|}e^{\frac{2\pi\mathrm{i}k(l+z)}{b}-\frac{\pi\mathrm{i}(a(l+z)-by)}{b}}
\overline{E}_{j-1}\biggl(\frac{a(l+z)}{b}-y\biggl)\nonumber\\
&&\quad+2a^{m}m!n!(a,b)\mathrm{sgn}(a)\sum_{j=1}^{n+1}\binom{m+n+1-j}{m}\frac{(-1)^{j}b^{n+1-j}}{(j-1)!\bigl(2\pi\mathrm{i}(k-\frac{a+b}{2})\bigl)^{m+n+2-j}}\nonumber\\
&&\qquad\times\sum_{l=1}^{|\frac{a}{(a,b)}|}e^{\frac{2\pi\mathrm{i}k(l+y)}{a}-\frac{\pi\mathrm{i}(b(l+y)-az)}{a}}
\overline{E}_{j-1}\biggl(\frac{b(l+y)}{a}-z\biggl).
\end{eqnarray}
By inserting \eqref{eq4.22}, \eqref{eq4.24} and \eqref{eq4.26} into \eqref{eq4.19}, we get that
\begin{eqnarray*}
&&e^{\pi\mathrm{i}(ax+bx+y+z)}\overline{E}_{m}(ax+y)\overline{E}_{n}(bx+z)+\delta_{0,m}\delta_{0,n}\mathrm{sgn}(ab)\delta_{\mathbb{Z}}(ax+y)\delta_{\mathbb{Z}}(bx+z)\nonumber\\
&&=\sideset{}{'}\sum_{k=-\infty}^{+\infty}\biggl(2b^{n}m!n!(a,b)\mathrm{sgn}(b)
\sum_{j=1}^{m+1}\binom{m+n+1-j}{n}\nonumber\\
&&\qquad\times\frac{(-1)^{j}a^{m+1-j}}{(j-1)!\bigl(2\pi\mathrm{i}((a,b)k-\frac{a+b}{2})\bigl)^{m+n+2-j}}\nonumber\\
&&\qquad\times\sum_{l=1}^{|\frac{b}{(a,b)}|}e^{\frac{2\pi\mathrm{i}(a,b)k(l+z)}{b}-\frac{\pi\mathrm{i}((a(l+z)-by)}{b}}
\overline{E}_{j-1}\biggl(\frac{a(l+z)}{b}-y\biggl)\nonumber\\
&&\quad+2a^{m}m!n!(a,b)\mathrm{sgn}(a)\sum_{j=1}^{n+1}\binom{m+n+1-j}{m}\nonumber\\
&&\qquad\times\frac{(-1)^{j}b^{n+1-j}}{(j-1)!\bigl(2\pi\mathrm{i}((a,b)k-\frac{a+b}{2})\bigl)^{m+n+2-j}}\nonumber\\
&&\qquad\times\sum_{l=1}^{|\frac{a}{(a,b)}|}e^{\frac{2\pi\mathrm{i}(a,b)k(l+y)}{a}-\frac{\pi\mathrm{i}(b(l+y)-az)}{a}}
\overline{E}_{j-1}\biggl(\frac{b(l+y)}{a}-z\biggl)\biggl)e^{2\pi\mathrm{i}(a,b)kx}\nonumber\\
&&\quad+\frac{(-1)^{n+1}2m!n!(a,b)^{m+n+2}e^{\pi\mathrm{i}(ax+bx+y+z)}}{a^{n+1}b^{m+1}(m+n+1)!}\overline{E}_{m+n+1}\biggl(\frac{by-az}{(a,b)}\biggl),
\end{eqnarray*}
from which and $2(a,b)\mid (a+b)$ we have
\begin{eqnarray}\label{eq4.27}
&&\overline{E}_{m}(ax+y)\overline{E}_{n}(bx+z)+\delta_{0,m}\delta_{0,n}\mathrm{sgn}(ab)(-1)^{ax+bx+y+z}\delta_{\mathbb{Z}}(ax+y)\delta_{\mathbb{Z}}(bx+z)\nonumber\\
&&=\sideset{}{'}\sum_{k=-\infty}^{+\infty}\biggl(2b^{n}m!n!\mathrm{sgn}(b)
\sum_{j=1}^{m+1}\binom{m+n+1-j}{n}\nonumber\\
&&\qquad\times\frac{(-1)^{j}a^{m+1-j}}{(j-1)!(a,b)^{m+n+1-j}(2\pi\mathrm{i}k)^{m+n+2-j}}\nonumber\\
&&\qquad\times\sum_{l=1}^{|\frac{b}{(a,b)}|}(-1)^{l}e^{\frac{2\pi\mathrm{i}(a,b)k(l+z)}{b}}
\overline{E}_{j-1}\biggl(\frac{a(l+z)}{b}-y\biggl)\nonumber\\
&&\quad+2a^{m}m!n!\mathrm{sgn}(a)\sum_{j=1}^{n+1}\binom{m+n+1-j}{m}\nonumber\\
&&\qquad\times\frac{(-1)^{j}b^{n+1-j}}{(j-1)!(a,b)^{m+n+1-j}(2\pi\mathrm{i}k)^{m+n+2-j}}\nonumber\\
&&\qquad\times\sum_{l=1}^{|\frac{a}{(a,b)}|}(-1)^{l}e^{\frac{2\pi\mathrm{i}(a,b)k(l+y)}{a}}
\overline{E}_{j-1}\biggl(\frac{b(l+y)}{a}-z\biggl)\biggl)e^{2\pi\mathrm{i}(a,b)kx}\nonumber\\
&&\quad+\frac{(-1)^{n+1}2m!n!(a,b)^{m+n+2}}{a^{n+1}b^{m+1}(m+n+1)!}\overline{E}_{m+n+1}\biggl(\frac{by-az}{(a,b)}\biggl).
\end{eqnarray}
Hence, we see from \eqref{eq2.1} and \eqref{eq4.27} that
\begin{eqnarray}\label{eq4.28}
&&\overline{E}_{m}(ax+y)\overline{E}_{n}(bx+z)+\delta_{0,m}\delta_{0,n}\mathrm{sgn}(ab)(-1)^{ax+bx+y+z}\delta_{\mathbb{Z}}(ax+y)\delta_{\mathbb{Z}}(bx+z)\nonumber\\
&&=-2b^{n}\mathrm{sgn}(b)
\sum_{j=1}^{m+1}\binom{m}{j-1}\frac{(-1)^{j}a^{m+1-j}}{(a,b)^{m+n+1-j}(m+n+2-j)}\nonumber\\
&&\qquad\times\sum_{l=1}^{|\frac{b}{(a,b)}|}(-1)^{l}\overline{E}_{j-1}\biggl(\frac{a(l+z)}{b}-y\biggl)
\overline{B}_{m+n+2-j}\biggl((a,b)x+\frac{(a,b)(l+z)}{b}\biggl)\nonumber\\
&&\quad-2a^{m}\mathrm{sgn}(a)\sum_{j=1}^{n+1}\binom{n}{j-1}\frac{(-1)^{j}b^{n+1-j}}{(a,b)^{m+n+1-j}(m+n+2-j)}\nonumber\\
&&\qquad\times\sum_{l=1}^{|\frac{a}{(a,b)}|}(-1)^{l}\overline{E}_{j-1}\biggl(\frac{b(l+y)}{a}-z\biggl)
\overline{B}_{m+n+2-j}\biggl((a,b)x+\frac{(a,b)(l+y)}{a}\biggl)\nonumber\\
&&\quad+\frac{(-1)^{n+1}2m!n!(a,b)^{m+n+2}}{a^{n+1}b^{m+1}(m+n+1)!}\overline{E}_{m+n+1}\biggl(\frac{by-az}{(a,b)}\biggl).
\end{eqnarray}
It is easily shown from \eqref{eq1.12}, \eqref{eq1.13}, \eqref{eq2.5}, the property of residue systems and the division algorithm that for $2\nmid a$ and $2\nmid b$,
\begin{eqnarray}\label{eq4.29}
&&\frac{1}{(a,b)^{m+n+1-j}}\sum_{l=1}^{|\frac{b}{(a,b)}|}(-1)^{l}
\overline{E}_{j-1}\biggl(\frac{a(l+z)}{b}-y\biggl)\nonumber\\
&&\qquad\times\overline{B}_{m+n+2-j}\biggl((a,b)x+\frac{(a,b)(l+z)}{b}\biggl)\nonumber\\
&&=\sum_{k=0}^{(a,b)-1}\sum_{l=1}^{|\frac{b}{(a,b)}|}(-1)^{k|\frac{b}{(a,b)}|+l}\overline{E}_{j-1}\biggl(\frac{a\bigl(k|\frac{b}{(a,b)}|+l+z\bigl)}{b}-y\biggl)\nonumber\\
&&\qquad\times\overline{B}_{m+n+2-j}\biggl(x+\frac{l+z+k|\frac{b}{(a,b)}|}{b}\biggl)\nonumber\\
&&=\sum_{l=1}^{|b|}(-1)^{l}\overline{E}_{j-1}\biggl(\frac{a(l+z)}{b}-y\biggl)\overline{B}_{m+n+2-j}\biggl(x+\frac{l+z}{b}\biggl).
\end{eqnarray}
Therefore, by applying \eqref{eq4.29} to the right-hand side of \eqref{eq4.28}, we get \eqref{eq2.20} and prove Theorem \ref{thm2.8}.

\noindent{\em{The proof of Theorem \ref{thm2.10}.}} Since $2\nmid(a+b)$, by inserting \eqref{eq4.22} and \eqref{eq4.26} into \eqref{eq4.19}, we see that
\begin{eqnarray*}
&&e^{\pi\mathrm{i}(ax+bx+y+z)}\overline{E}_{m}(ax+y)\overline{E}_{n}(bx+z)+\delta_{0,m}\delta_{0,n}\mathrm{sgn}(ab)\delta_{\mathbb{Z}}(ax+y)\delta_{\mathbb{Z}}(bx+z)\nonumber\\
&&=\sum_{k=-\infty}^{+\infty}\biggl(2b^{n}m!n!(a,b)\mathrm{sgn}(b)
\sum_{j=1}^{m+1}\binom{m+n+1-j}{n}\nonumber\\
&&\qquad\times\frac{(-1)^{j}a^{m+1-j}}{(j-1)!\bigl(2\pi\mathrm{i}((a,b)k-\frac{a+b}{2})\bigl)^{m+n+2-j}}\nonumber\\
&&\qquad\times\sum_{l=1}^{|\frac{b}{(a,b)}|}e^{\frac{2\pi\mathrm{i}(a,b)k(l+z)}{b}-\frac{\pi\mathrm{i}((a(l+z)-by)}{b}}
\overline{E}_{j-1}\biggl(\frac{a(l+z)}{b}-y\biggl)\nonumber\\
&&\quad+2a^{m}m!n!(a,b)\mathrm{sgn}(a)\sum_{j=1}^{n+1}\binom{m+n+1-j}{m}\nonumber\\
&&\qquad\times\frac{(-1)^{j}b^{n+1-j}}{(j-1)!\bigl(2\pi\mathrm{i}((a,b)k-\frac{a+b}{2})\bigl)^{m+n+2-j}}\nonumber\\
&&\qquad\times\sum_{l=1}^{|\frac{a}{(a,b)}|}e^{\frac{2\pi\mathrm{i}(a,b)k(l+y)}{a}-\frac{\pi\mathrm{i}(b(l+y)-az)}{a}}
\overline{E}_{j-1}\biggl(\frac{b(l+y)}{a}-z\biggl)\biggl)e^{2\pi\mathrm{i}(a,b)kx}.
\end{eqnarray*}
It follows that
\begin{eqnarray}\label{eq4.30}
&&\overline{E}_{m}(ax+y)\overline{E}_{n}(bx+z)+\delta_{0,m}\delta_{0,n}\mathrm{sgn}(ab)(-1)^{ax+bx+y+z}\delta_{\mathbb{Z}}(ax+y)\delta_{\mathbb{Z}}(bx+z)\nonumber\\
&&=\sum_{k=-\infty}^{+\infty}\biggl(2b^{n}m!n!\mathrm{sgn}(b)
\sum_{j=1}^{m+1}\binom{m+n+1-j}{n}\nonumber\\
&&\qquad\times\frac{(-1)^{j}a^{m+1-j}}{(j-1)!(a,b)^{m+n+1-j}\bigl(2\pi\mathrm{i}(k-\frac{1}{2})\bigl)^{m+n+2-j}}\nonumber\\
&&\qquad\times\sum_{l=1}^{|\frac{b}{(a,b)}|}(-1)^{l}e^{\frac{2\pi\mathrm{i}(a,b)(k-\frac{1}{2})(l+z)}{b}}
\overline{E}_{j-1}\biggl(\frac{a(l+z)}{b}-y\biggl)\nonumber\\
&&\quad+2a^{m}m!n!\mathrm{sgn}(a)\sum_{j=1}^{n+1}\binom{m+n+1-j}{m}\nonumber\\
&&\qquad\times\frac{(-1)^{j}b^{n+1-j}}{(j-1)!(a,b)^{m+n+1-j}\bigl(2\pi\mathrm{i}(k-\frac{1}{2})\bigl)^{m+n+2-j}}\nonumber\\
&&\qquad\times\sum_{l=1}^{|\frac{a}{(a,b)}|}(-1)^{l}e^{\frac{2\pi\mathrm{i}(a,b)(k-\frac{1}{2})(l+y)}{a}}
\overline{E}_{j-1}\biggl(\frac{b(l+y)}{a}-z\biggl)\biggl)\nonumber\\
&&\qquad\times e^{2\pi\mathrm{i}(a,b)(k-\frac{1}{2})x}.
\end{eqnarray}
Applying \eqref{eq2.2} to the right-hand side of \eqref{eq4.30}, we discover that
\begin{eqnarray}\label{eq4.31}
&&\overline{E}_{m}(ax+y)\overline{E}_{n}(bx+z)+\delta_{0,m}\delta_{0,n}\mathrm{sgn}(ab)(-1)^{ax+bx+y+z}\delta_{\mathbb{Z}}(ax+y)\delta_{\mathbb{Z}}(bx+z)\nonumber\\
&&=b^{n}\mathrm{sgn}(b)
\sum_{j=1}^{m+1}\binom{m}{j-1}\frac{(-1)^{j}a^{m+1-j}}{(a,b)^{m+n+1-j}}\sum_{l=1}^{|\frac{b}{(a,b)}|}(-1)^{l}
\overline{E}_{j-1}\biggl(\frac{a(l+z)}{b}-y\biggl)\nonumber\\
&&\qquad\times\overline{E}_{m+n+1-j}\biggl((a,b)x+\frac{(a,b)(l+z)}{b}\biggl)\nonumber\\
&&\quad+a^{m}\mathrm{sgn}(a)\sum_{j=1}^{n+1}\binom{n}{j-1}\frac{(-1)^{j}b^{n+1-j}}{(a,b)^{m+n+1-j}}
\sum_{l=1}^{|\frac{a}{(a,b)}|}(-1)^{l}\overline{E}_{j-1}\biggl(\frac{b(l+y)}{a}-z\biggl)\nonumber\\
&&\qquad\times\overline{E}_{m+n+1-j}\biggl((a,b)x+\frac{(a,b)(l+y)}{a}\biggl).
\end{eqnarray}
Note that from \eqref{eq1.13}, \eqref{eq4.16} and the
division algorithm, we see for $2\nmid (a+b)$ that
\begin{eqnarray}\label{eq4.32}
&&\frac{1}{(a,b)^{m+n+1-j}}\sum_{l=1}^{|\frac{b}{(a,b)}|}(-1)^{l}
\overline{E}_{j-1}\biggl(\frac{a(l+z)}{b}-y\biggl)\nonumber\\
&&\qquad\times\overline{E}_{m+n+1-j}\biggl((a,b)x+\frac{(a,b)(l+z)}{b}\biggl)\nonumber\\
&&=\sum_{k=0}^{(a,b)-1}\sum_{l=1}^{|\frac{b}{(a,b)}|}(-1)^{k|\frac{b}{(a,b)}|+l}\overline{E}_{j-1}\biggl(\frac{a\bigl(k|\frac{b}{(a,b)}|+l+z\bigl)}{b}-y\biggl)\nonumber\\
&&\qquad\times\overline{E}_{m+n+1-j}\biggl(x+\frac{l+z+k|\frac{b}{(a,b)}|}{b}\biggl)\nonumber\\
&&=\sum_{l=1}^{|b|}(-1)^{l}\overline{E}_{j-1}\biggl(\frac{a(l+z)}{b}-y\biggl)\overline{E}_{m+n+1-j}\biggl(x+\frac{l+z}{b}\biggl).
\end{eqnarray}
Thus, by applying \eqref{eq4.32} to the right-hand side of \eqref{eq4.31}, we achieve \eqref{eq2.22} and conclude the proof of Theorem \ref{thm2.10}.

\end{document}